\documentstyle[11pt]{article}

\title{ORDINAL DISTANCES IN TRANSFINITE GRAPHS}
\author{A. H. Zemanian}
\date{}

\begin{document}
\maketitle
\baselineskip21pt

{\ Abstract --- An ordinal-valued 
metric, taking its values in the set $\aleph_{1}$ of all countable ordinals,
can be assigned to a metrizable set $\cal M$
of nodes in any transfinite graph.  $\cal M$ contains all the 
nonsingleton nodes, as well as certain singleton nodes.
Moreover, this yields a 
graphical realization of Cantor's countable ordinals, as well as
of the Aristotelian ideas of ``potential'' and ``actual'' infinities,
the former being represented by the arrow ranks and the latter 
by the ordinal ranks of transfiniteness.  This construct also extends 
transfinitely the ideas of nodal eccentricities, radii, 
diameters, centers, peripheries, and blocks for graphs,
and the following generalizations are established.
With $\nu$ denoting the rank of 
a $\nu$-graph ${\cal G}^{\nu}$, the $\nu$-nodes of ${\cal G}^{\nu}$
comprise the center of a larger $\nu$-graph.  Also, when there are only 
finitely many $\nu$-nodes and when those $\nu$-nodes are ``pristine''
in the sense that they do not embrace nodes of lower ranks, the infinitely 
many nodes of all ranks have eccentricities of the form 
$\omega^{\nu}\cdot p$, where $\omega$ is the first transfinite ordinal and $p$
lies in a finite set of natural numbers.  Furthermore,
the center is contained in a single block of highest rank.  Also, when 
each loop of the $\nu$-graph is confined within a $(\nu-1)$-section,
the center either is a single node of highest rank, or is the set of 
internal nodes of a $(\nu-1)$-section, or is
the union of the latter two kinds of centers.\\

} 

\vspace{.07in}
1. Introduction

The idea of distances in connected finite graphs has been quite fruitful,
with much research directed toward both theory and applications.
See, for example, \cite{bu}, \cite{b-h}, 
\cite{c-l}, 
and the references therein.  Such distances are given by a metric that 
assigns to each pair of nodes the minimum number of branches among all paths 
connecting those two nodes.  Thus, the metric takes its values in the set 
$\aleph_{0}$ of natural numbers.  That 
distance considerations can be so fruitful in the theory of 
finite graphs inspires the question of whether distance constructs 
can be devised for transfinite graphs.  Transfinite graphs 
\cite{tgen}, \cite{pg} represent a 
generalization of graphs that is roughly analogous to Cantor's extension 
of the natural numbers to the transfinite ordinals.

An affirmative answer to that question was achieved in \cite{c-z}, wherein 
a real-valued metric was devised for the purpose of ascertaining
limit points at infinite extremities of a conventionally infinite,  
electrical, resistive network, through which points
electrical current could flow into other such networks.
This construct was extended to higher ranks of transfiniteness 
\cite{pg} with an infinite hierarchy of metrics, one for each rank of 
transfiniteness.  These metrics take their values in the nonnegative 
real line, are quite different from
the standard branch-count metric mentioned above, require a variety of 
restriction such as local finiteness, and do not reduce to the 
branch-count metric for finite graphs.

Let us emphasize why, we feel, 
that for this paper it is inappropriate 
to use a real-valued metric
that makes infinite extremities of a conventionally infinite graph 
look as though they are only finitely distant from any node of the graph.
If branch counts are to determine distances between nodes, then no node 
is closer to any infinite extremity than any other node.  Thus,
all nodes in a conventionally infinite graph
should be viewed as equally distant from any extremity, and that 
distance should be $\omega$, the first transfinite ordinal.  This 
property cannot be avoided if branch counts are to prevail.  

The problem attacked in this work is the invention of a single metric 
that extends the standard branch-count metric to transfinite graphs,
one that holds for all ranks of transfiniteness, and reduces to
the standard branch-count metric for finite graphs.
In closer analogy to Cantor's work, the metric proposed in this paper assigns 
countable ordinals to pairs of nodes in a connected transfinite graph;
that is, it takes its values in the set $\aleph_{1}$ of all countable
ordinals.  Moreover, the metric is applicable even when the graph  
is not locally finite and may even have uncountably many branches.

As a consequence, the ideas of nodal eccentricities, 
radii, diameters, centers, peripheries, and blocks 
are herein extended to transfinite graphs. However, to do so, the 
set $\aleph_{1}$ has to be enlarged by inserting an ``arrow rank''
\cite[page 4]{tgen}, \cite[page 4]{pg} immediately preceding each limit 
ordinal.  These arrow ranks reflect the Aristotelian idea of a 
``potential infinity'' as distinct from the other Aristotelian 
idea of an ``actual infinity'' typified by the ordinals.

Several theorems concerning these ideas are proven, some of which 
lift results concerning finite graphs to transfinite graphs.
For example, with ${\cal G}^{\nu}$ henceforth denoting a transfinite graph of 
rank $\nu$ (i.e., a $\nu$-graph) \cite[Chap. 2]{tgen},
the transfinite radius and diameter of ${\cal G}^{\nu}$
are related, but now in a more complicated way (Theorem 7.2).
Nodes of highest rank (i.e., $\nu$-nodes) are shown to comprise the 
center of a larger $\nu$-graph (Corollary 7.5).  A finite range of 
possible transfinite ordinals for the eccentricities of the
nodes of ${\cal G}^{\nu}$ is established (Corollary 8.10).
The center of ${\cal G}^{\nu}$ lies in a block of highest rank
(Theorem 9.5), and that center is identified either as a single 
$\nu$-node or as a certain infinite set of nodes of all ranks
(Theorem 10.3).  To obtain these and other results, we employ some 
restrictions on ${\cal G}^{\nu}$, which are introduced when needed and 
are then assumed to hold throughout the rest of the paper.

Various properties of transfinite graphs are used in this work.
These can be found in the book \cite{tgen}.  A simplified but more 
restrictive rendition of the subject is given in \cite{pg}.  We will work 
in the generality of \cite{tgen} and will refer to specific pages in that 
book when invoking various concepts and results. In this paper,
we do not allow any branch to be a self-loop;  
thus, every branch is incident to two different nodes.
However, parallel branches are allowed.
We use the standard notations for ordinals and cardinals \cite{ab}.

Furthermore, any transfinite node $x^{\alpha}$ may (but need not)
contain exactly one node of lower rank $x^{\beta}$ $(\beta <\alpha)$;
$x^{\beta}$ in turn may contain exactly one 
other node $x^{\gamma}$ $(\gamma <\beta)$, and so forth 
through finitely many decreasing ranks.  We say that $x^{\alpha}$
{\em embraces} itself and $x^{\beta}$, $x^{\gamma}$, $\ldots\;$, as well.
On the other hand, if $x^{\alpha}$ is not embraced by a node of
higher rank, we call $x^{\alpha}$ a {\em maximal} node.  It is the 
maximal nodes we will be primarily 
concerned with because connectedness to $x^{\alpha}$ implies 
connectedness to $x^{\beta}$, $x^{\gamma}$, $\dots\;$, as well. 
Rather than repeating the adjective ``maximal,''
we let it be understood throughout 
that any node discussed is maximal unless the 
opposite is explicitly stated.
This implies that different (maximal) nodes must be ``totally disjoint,''
that is, they embrace no common elements \cite[Lemma 2.2-1]{tgen}.

Throughout this work we restrict the rank $\nu$ of ${\cal G}^{\nu}$
to $1\leq\nu\leq\omega$, $\nu\neq\vec{\omega}$.  (Here, 
$\vec{\omega}$ is the arrow rank immediately preceding $\omega$
\cite[Sec. 2.3]{tgen}.)  All our results can be extended to many 
ordinal ranks higher than $\omega$ by modifying our arguments in obvious ways.
However, it is not apparent whether this extension can be made 
throughout all the ordinal ranks in $\aleph_{1}$ \cite[Sec. 2.5]{tgen}.

Given any ordinal rank $\rho$ no higher than $\nu$, two nodes 
(resp. two branches) of ${\cal G}^{\nu}$ are said to be $\rho$-connected
if there is a path of rank $\rho$ or less that terminates at those two
nodes (resp. terminates at 0-nodes of those branches).\footnote{Paths
of various ranks are defined in the next section.  See also 
\cite[Chap. 2]{tgen}.} If this holds for $\rho=\nu$ and for all 
branches in ${\cal G}^{\nu}$, ${\cal G}^{\nu}$ is said to be 
$\nu$-connected.  We always assume that ${\cal G}^{\nu}$ is $\nu$-connected.

\vspace{.07in}
2. Lengths of Paths

Throughout this paper, we use the natural sum of transfinite ordinals 
to obtain the normal expansion of that sum \cite[pages 354-355]{ab}.

\vspace{.1in}
\noindent
{\em 0-Paths:}

A (nontrivial) 0-path $P^{0}$ is an alternating sequence 
\begin{equation}
P^{0}\;=\;\{\ldots,x_{m}^{0},b_{m},x_{m+1}^{0},b_{m+1},\ldots\} \label{2.1}
\end{equation}
of branches $b_{m}$ and conventional nodes $x_{m}^{0}$ (also called
``0-nodes'') in which no term repeats and each branch is incident
to the two 0-nodes adjacent to it in the sequence. If the sequence
terminates on either side, it terminates at a 0-node.
This is the conventional 
definition of a path.  (The 0-nodes of (\ref{2.1}) need not be 
maximal when $P^{0}$ occurs within a transfinite graph.)   
When $P^{0}$ is one-ended (i.e., one-way infinite), its length is
defined to be $|P^{0}|={\omega}$.  When $P^{0}$ is endless 
(i.e., two-way infinite),
its length is taken to be $|P^{0}|={\omega}\cdot 2$.  If $P^{0}$ is two-ended
(i.e., has only finitely many 0-nodes), we set 
$|P^{0}|=\tau_{0}$, where $\tau_{0}$ is the number of 
branches in $P^{0}$.  We might motivate these definitions by 
noting that we are using $\omega$ to denote the infinity of branches 
in a one-ended 0-path and using $\omega\cdot 2$ to represent to fact 
that an endless 0-path is the union of two one-ended paths.
Equivalently, we can identify 
$\omega$ with each 0-tip traversed;  a one-ended 
0-path has one 0-tip, and an endless 0-path has two 0-tips---hence, the 
length $\omega\cdot 2$. (See \cite[page 20]{tgen} for 
the definition of a 0-tip.)

\vspace{.1in}
\noindent
{\em 1-Paths.}

A (nontrivial) 1-path $P^{1}$ \cite[page 28]{tgen} is an alternating sequence 
\begin{equation}
P^{1}\;=\;\{\ldots,x_{m}^{1},P_{m}^{0},x_{m+1}^{1},P_{m+1}^{0},\ldots\}  \label{2.2}
\end{equation}
of 1-nodes $x_{m}^{1}$ and 0-paths $P_{m}^{0}$ that represents 
a tracing through a transfinite graph of rank 1 or greater in which no 
node is met more than once in the tracing. 
If the sequence terminates on either side,
it terminates at a 0-node or 1-node.
See \cite[page 28]{tgen} for the full definition of a 1-path.
The length $|P^{1}|$ of $P^{1}$ is defined as follows.  When $P^{1}$
is one-ended, $|P^{1}|=\omega^{2}$, and, when $P^{1}$ is endless, $|P^{1}|=
\omega^{2}\cdot 2$.  When $P^{1}$ is two-ended (i.e., when it has 
only finitely many 1-nodes), we set $|P^{1}|=\sum_{m}|P_{m}^{0}|$, where 
the sum is over the finitely many 0-paths $P_{m}^{0}$ in (\ref{2.2});
thus, in this case, $|P^{1}|=\omega\cdot \tau_{1}+\tau_{0}$, where
$\tau_{1}$ is the number of 0-tips $P^{1}$ traverses, and $\tau_{0}$
is the number of branches in all the 0-paths in (\ref{2.2}) 
that are two-ended. 
It is important here to write $|P^{1}|$ as indicated and not as $\tau_{0}+
\omega\cdot \tau_{1}$ because ordinal addition is not commutative
\cite[page 327]{ab}.  Thus, $\omega\cdot \tau_{1}+\tau_{0}$ takes into 
account the lengths of all the 0-paths in (\ref{2.2}), but 
$\tau_{0}+\omega\cdot\tau_{1}$ fails to do so.  As specified above, 
$\omega\cdot \tau_{1}+\tau_{0}$ is the ``normal expansion''
\cite[pages 354-355]{ab} of $|P^{1}|$.
     
\vspace{.1in}
\noindent
{\em $\mu$-Paths:}

Now, let $\mu$ be any positive natural number. A $\mu$-path 
\cite[page 33]{tgen} is an 
alternating sequence
\begin{equation}
P^{\mu}\;=\;\{\ldots,x_{m}^{\mu},P_{m}^{\alpha_{m}},x_{m+1}^{\mu},
P_{m+1}^{\alpha_{m+1}},\ldots\}   \label{2.3}
\end{equation}
of $\mu$-nodes $x_{m}^{\mu}$ and $\alpha_{m}$-paths $P_{m}^{\alpha_{m}}$,
where $0\leq \alpha_{m}<\mu$. (The natural numbers 
$\alpha_{m}$ may vary with $m$, and the 
$\mu$-nodes need not be maximal.)  As before, $P^{\mu}$ represents 
a tracing through a transfinite graph of rank $\mu$ or larger 
in which no node is met more than once in the tracing.
Termination on either side of (\ref{2.3}) occurs at a node of rank $\mu$ or less. 
When $P^{\mu}$ is one-ended, its length $|P^{\mu}|$ is defined to be 
$\omega^{\mu+1}$, and, when $P^{\mu}$ is endless , we set
$|P^{\mu}|=\omega^{\mu+1}\cdot 2$.  
When, however, $P^{\mu}$ is two-ended (i.e., has only finitely many
$\mu$-nodes), we set $|P^{\mu}|=\sum_{m}|P_{m}^{\alpha_{m}}|$,
where as always this sum denotes a normal expansion of an ordinal
obtained through a natural summation of ordinals 
\cite[pages 354-355]{ab}. Recursively, this gives
\begin{equation}
|P^{\mu}|\;=\;\omega^{\mu}\cdot \tau_{\mu}+\omega^{\mu-1}\cdot \tau_{\mu-1}+\ldots+\omega\cdot \tau_{1}+\tau_{0}, \label{2.4}
\end{equation}
where $\tau_{\mu},\tau_{\mu-1},\ldots,\tau_{0}$ are natural numbers.
$\tau_{\mu}$ is the number of $(\mu-1)$-tips among all 
the one-ended and endless $(\mu-1)$-paths (i.e., when $\alpha_{m}=\mu-1$)
appearing in (\ref{2.3});  $\tau_{\mu}$ is not 0.
For $k=\mu-1, \mu-2,\ldots,1$, we set $\tau_{k}$ 
equal to the number of $k-1$-tips generated by these recursive definitions.
Finally, $\tau_{0}$ is one-half 
the number of elementary tips \cite[page 9]{tgen}
generated recursively by these definitions.  
Thus, $\tau_{0}$ is a number of branches because each branch has exactly two elementary tips. 
Any $\tau_{k}$ $(k<\mu$) can be 0.

{\bf Example 2.1.}  Let $P^{3}$ be the two-ended 3-path:
\[ P^{3}\;=\;\{x_{1}^{2},P^{2}_{1},x_{2}^{3},P_{2}^{2},x_{3}^{3},
P_{3}^{2},x_{4}^{3}\} \]
Here, $P_{1}^{2}$ is assumed to be a one-ended 2-path terminating on the left 
with $x_{1}^{2}$ and 
reaching $x_{2}^{3}$
through a 2-tip.  Hence, $|P_{1}^{2}|=\omega^{3}$.  We take $P_{2}^{2}$
to be the two-ended 2-path 
\[ P_{2}^{2}\;=\;\{y_{1}^{2},Q_{1}^{1},y_{2}^{2},Q_{2}^{0},y_{3}^{2}\}, \]
where $y_{1}^{2}$ and $y_{3}^{2}$ are members of $x_{2}^{3}$ 
and $x_{3}^{3}$ respectively. $Q_{1}^{1}$ is an endless 1-path reaching the 
2-nodes $y_{1}^{2}$ and $y_{2}^{2}$ with 1-tips, and $Q_{2}^{0}$
is a finite 0-path with four branches, whose terminal 0-nodes are
members of $y_{2}^{2}$ and $y_{3}^{2}$.  Hence,
$|Q_{1}^{1}|=\omega^{2}\cdot 2$ and $|Q_{2}^{0}|=4$.  Finally, 
we take $P_{3}^{2}$ to be an endless 2-path reaching $x_{3}^{3}$
and $x_{4}^{3}$ through 2-tips.  Hence, $|P_{3}^{2}|=\omega^{3}\cdot 2$.

Altogether then, with a rearrangement of the following ordinal sum 
to get a normal-expansion, we may write
\begin{eqnarray*}
|P^{3}|\;&=&\;|P_{1}^{2}|+|P_{2}^{2}|+|P_{3}^{2}|\\
         &=&\;\omega^{3}+|Q_{1}^{1}|+|Q_{2}^{0}|+\omega^{3}\cdot 2\\
         &=&\;\omega^{3}+\omega^{2}\cdot 2+4+\omega^{3}\cdot 2\\
         &=&\;\omega^{3}\cdot 3 +\omega^{2}\cdot 2 + 4
\end{eqnarray*}
$\Box$

\vspace{.1in}
\noindent
{\em $\vec{\omega}$-paths:}

$\vec{\omega}$-paths occur within paths of ranks $\omega$ and higher, but they are 
never two-ended \cite[pages 40-41]{tgen}.  The length of an $\vec{\omega}$-path
$P^{\vec{\omega}}$ is defined to be $|P^{\vec{\omega}}|=\omega^{\omega}$ 
when $P^{\vec{\omega}}$
is one-ended, and $|P^{\vec{\omega}}|=\omega^{\omega}\cdot 2$ 
when $P^{\vec{\omega}}$ is endless.

\vspace{.1in}
\noindent
{\em $\omega$-paths:}

A (nontrivial) $\omega$-path $P^{\omega}$ \cite[page 44]{tgen}
\begin{equation}
P^{\omega}\;=\;\{\ldots, x_{m}^{\omega},P_{m}^{\alpha_{m}},x_{m+1}^{\omega},
P_{m+1}^{\alpha_{m+1}},\ldots\}  \label{2.5}
\end{equation}
is an alternating sequence of (not necessarily maximal) $\omega$-nodes 
$x_{m}^{\omega}$ and $\alpha_{m}$-paths $P_{m}^{\alpha_{m}}$ 
$(0\leq\alpha_{m}\leq\vec{\omega} )$ that represents a 
tracing through a graph of rank $\omega$ (or larger) in 
which no node is met more than once 
and a termination on either side is at a node of rank 
$\omega$ or less.  By definition, 
when $P^{\omega}$ is one-ended, $|P^{\omega}|=\omega^{\omega+1}$;
also, when $P^{\omega}$ is endless, $|P^{\omega}|= \omega^{\omega+1}\cdot 2$.
When $P^{\omega}$ is two-ended (i.e., has only finitely many $\omega$-nodes),
we set 
\begin{equation}
|P^{\omega}|\;=\;\sum_{m}|P_{m}^{\alpha_{m}}|\;=\;\omega^{\omega}\cdot \tau_{\omega}
+\sum_{k=0}^{\infty}\omega^{k}\cdot \tau_{k}  \label{2.6}
\end{equation}
with the natural summation being understood.
Here, $\tau_{\omega}$ is the number of 
$\vec{\omega}$-tips among all the one-ended 
and endless $\vec{\omega}$-paths appearing as elements 
of $P_{m}^{\alpha_{m}}$ in 
(\ref{2.5}) (i.e., when $\alpha_{m}=\vec{\omega}$); $\tau_{\omega}$ is not 0. 
On the other hand, 
the $\tau_{k}$
are determined recursively, as they are in (\ref{2.4}).  There are 
only finitely many nonzero terms in the summation within 
(\ref{2.6}) because there are 
only finitely many paths $P_{m}^{\alpha_{m}}$ in a two-ended $\omega$-path
and each $|P_{m}^{\alpha_{m}}|$ is a finite sum as in (\ref{2.4}). 

An immediate result of all these definitions is the following.

{\bf Lemma 2.2.}  {\em If $Q^{\beta}$ is a subpath of a $\gamma$-path
$P^{\gamma}$ $(0\leq \beta\leq\gamma)$, then $|Q^{\beta}|\leq|P^{\gamma}|$.}

It is easy to add ordinals when they are in 
normal-expansion form---simply
add their corresponding coefficients.  Thus, the length of the union of 
two paths that are totally disjoint except for incidence at a terminal 
node (a ``series connection'') is obtained by adding their lengths 
in normal expansion form. 
Similarly, if $Q$ is a 
proper subpath of $P$, the part
$P\backslash Q$ of $P$ that is not in $Q$ has the total length $|P|-|Q|$, 
which is obtained by subtracting the coefficients of $|Q|$ from the 
corresponding coefficients of $|P|$. 

\vspace{.07in}
3. Metrizable Sets of Nodes

In a connected finite graph, for every two nodes there is at least 
one path terminating at them. This is not in general true for 
transfinite graphs. 

{\bf Example 3.1.} The 1-graph of Fig. 1 provides an example.  
In that graph, $x_{a}^{1}$ (resp. $x_{b}^{1}$) is a nonsingleton 1-node 
containing the 0-tip $t_{a}^{0}$ (resp. $t_{b}^{0}$) for the 
one-ended path of $a_{k}$ branches (resp. $b_{k}$ branches)
and also embracing an elementary tip of branch $d$ (resp. $e$).  There
are, in addition, uncountably many 0-tips for paths that alternate 
infinitely often between the $a_{k}$ and $b_{k}$ branches by passing through 
$c_{k}$ branches;  those tips are contained in singleton 1-nodes, one for each.
$x_{abc}^{1}$ denotes one such singleton 1-node;  the others are not shown.
Note that there is no path connecting $x_{abc}^{1}$ to $x_{a}^{1}$
(or to any other 1-node) because any tracing between $x_{abc}^{1}$
and $x_{a}^{1}$ must repeat 0-nodes.  Thus, our definition (given in the 
next section) of the distance
between two nodes as the minimum path length for all paths connecting 
those nodes cannot be applied to $x_{abc}^{1}$ and $x_{a}^{1}$.
We seek some means of applying this distance concept to 
at least some pairs of nodes. $\Box$

To this end, we impose the following condition on the transfinite 
graph ${\cal G}^{\nu}$, which is understood to hold henceforth.

{\bf Condition 3.2.}  {\em If two tips (perforce of ranks less than $\nu$
and possibly differing) are nondisconnectable,\footnote{Two tips
are called {\em nondisconnectable} if their representative 
(one-ended) paths continue to meet no matter how far along
the representative paths one proceeds \cite[page 58]{tgen}.  
Two tips are called {\em disconnectable}
if they have representative paths that are totally disjoint.}
then either they are shorted together (i.e., are embraced by the same node)
or at least one of them is open (i.e., is the sole member of a 
singleton node).}

The 1-graph of Fig. 1 satisfies this condition.

The following results ensue:  As specified above, 
${\cal G}^{\nu}$ is $\nu$-{\em connected}, which means that
for any two branches there is a two-ended path $P^{\rho}$
of some rank $\rho$ $(\rho\leq \nu)$ that meets those two branches.
Nevertheless, there may be two nodes not having any path 
that meets them (i.e., the two nodes are not $\nu$-connected).
For instance, the 1-graph of Fig. 1 is 1-connected, but there is no path that 
meets $x_{a}^{1}$ and $x_{abc}^{1}$.  Now, as will be
established by Lemma 3.3 below, if ${\cal G}^{\nu}$ satisfies
Condition 3.2, then, for any two nonsingleton nodes, there will be
at least one two-ended path terminating at them.  As a result, we will be able
to define distances between nonsingleton nodes.   
Furthermore, some singleton nodes may 
be amenable to such distance measurements, as well.  To test 
this, we need merely append a new branch $b$ to a singleton node 
$x^{\alpha}$ by adding an elementary tip of $b$ to $x^{\alpha}$
as an embraced elementary tip
to get a nonsingleton node $\hat{x}^{\alpha}$, with the other
elementary tip of $b$ left open (i.e., $b$ is added as an end branch)---and
then check to see if Condition 3.2 is maintained.  More generally,
with ${\cal G}^{\nu}$ being $\nu$-connected
and satisfying Condition 3.2, let ${\cal M}$ 
be a set consisting 
of all the nonsingleton (maximal) nodes in ${\cal G}^{\nu}$ 
and possibly other singleton 
(maximal) nodes having the property that, if end branches
are appended to those singleton nodes simultaneously, Condition 3.2 is 
still satisfied by the resulting network.  
Any such set 
${\cal M}$ will be call a {\em metrizable} set of nodes.

{\bf Lemma 3.3.}  {\em Let ${\cal M}$ be a 
metrizable set of nodes in ${\cal G}^{\nu}$.  Then, 
for any two nodes of ${\cal M}$, there exists a two-ended path terminating
to those nodes.}

{\bf Proof.}  Let $x_{a}^{\alpha}$ and $x_{b}^{\beta}$ be two different
nodes in ${\cal M}$.  Since they are maximal, they must be totally
disjoint.  Then, by Condition 3.2, any tip in $x_{a}^{\alpha}$
is disconnectable from every tip in $x_{b}^{\beta}$;  indeed, if they were 
nondisconnectable, they would have to be shorted together,
making $x_{a}^{\alpha}$ and $x_{b}^{\beta}$ the same node.  Thus,
we can choose a representative path $P_{a}$ for that tip in $x_{a}^{\alpha}$
that is totally disjoint from a representative path $P_{b}$ for a tip in 
$x_{b}^{\beta}$.  By the definition of $\nu$-connectedness, there
will be a path $P_{ab}$ connecting a branch of $P_{a}$ and a branch of
$P_{b}$.  By \cite[Corollary 3.5-4]{tgen}, there is in the subgraph 
$P_{a}\cup P_{ab}\cup P_{b}$ induced by the branches of those three
paths a two-ended path terminating at $x_{a}^{\alpha}$ and 
$x_{b}^{\beta}$.  $\Box$

{\bf Example 3.4.}  For an illustration, remove in Fig. 1 branches $d$ and $e$
along with the 0-nodes $y_{a}^{0}$ and $y_{b}^{0}$, thereby making
$x_{a}^{1}$ and $x_{b}^{1}$ singleton 1-nodes.  Then, the remaining 0-nodes 
along with $x_{a}^{1}$ and $x_{b}^{1}$ comprise a metrizable set.
Also, those 0-nodes along with $x_{abc}^{1}$ comprise another 
metrizable set.  However, those 0-nodes along with $x_{a}^{1}$,
$x_{b}^{1}$, and $x_{abc}^{1}$ do not comprise a metrizable set.
$\Box$

\vspace{.07in}
4. Distances Between Nodes

Our objective now is to define ordinal distances between nodes 
whereby the metric axioms are satisfied.  
Let ${\cal M}$ 
be a metrizable 
set of nodes in ${\cal G}^{\nu}$.  
We define the distance 
function $d\! : {\cal M}\times {\cal M}\leadsto \aleph_{1}$ as follows:
If $x_{a}^{\alpha}$ and $x_{b}^{\beta}$ are different nodes
in $\cal M$, we set 
\begin{equation}
d(x_{a}^{\alpha},x_{b}^{\beta})\;=\;\min\{|P_{ab}|\! : 
P_{ab}\;{\rm is \;a\; two\;ended \;path \;terminating \;at\; } 
x_{a}^{\alpha} \;{\rm and}\; x_{b}^{\beta}\}.  \label{4.1} 
\end{equation}
If $x_{a}^{\alpha}=x_{b}^{\beta}$, 
we set $d(x_{a}^{\alpha},x_{b}^{\beta})=0$. 
By our constructions in Sec. 2, $|P_{ab}|$ is a countable ordinal no larger
than $\omega^{\omega}\cdot k$, where $k$ is a natural number.
Moreover, any set of ordinals is well-ordered and thus has a least member.
Therefore, the minimum indicated in (\ref{4.1}) exists, and
is a countable ordinal.

Obviously, $d(x_{a}^{\alpha},x_{b}^{\beta})>0$ if 
$x_{a}^{\alpha}\not= x_{b}^{\beta}$.
Moreover, $d(x_{a}^{\alpha},x_{b}^{\beta})=d(x_{b}^{\beta},x_{a}^{\alpha})$.  
It remains to prove the triangle inequality; namely, if $x_{a}^{\alpha}$,
$x_{b}^{\beta}$, and $x_{c}^{\gamma}$ are any three (maximal) 
nodes in $\cal M$, then
\begin{equation}
d(x_{a}^{\alpha},x_{b}^{\beta})\;\leq\; d(x_{a}^{\alpha},x_{c}^{\gamma})\,+\, 
d(x_{c}^{\gamma},x_{b}^{\beta}).  \label{4.2}
\end{equation}

This is easily done by invoking Lemma 2.2 and using 
\cite[Corollary 3.5-4]{tgen}, whose rather long proof needs Condition 3.2.

{\bf Proposition 4.1.} {\em $d$ satisfies the metric axioms.}

Clearly, $d$ reduces to the standard (branch-count) distance function 
when ${\cal G}^{\nu}$ is replaced by a finite graph.
We have achieved one of the objectives of this paper by showing that the 
branch-count distance function can be extended transfinitely
to any metrizable set of nodes in ${\cal G}^{\nu}$.

{\bf Example. 4.2.}  For the 1-graph of Fig. 1 and with $\cal M$ 
consisting of all the 0-nodes along with $x_{a}^{1}$ and $x_{b}^{1}$,
we have $d(x_{1}^{0},x_{2}^{0})=1$, $d(x_{1}^{0},x_{a}^{1})
=d(x_{1}^{0},x_{b}^{1})=\omega$, $d(x_{1}^{0},y_{a}^{0})=
d(x_{1}^{0},y_{b}^{0})=\omega+1$,
and $d(y_{a}^{0},y_{b}^{0})=\omega\cdot 2+2$.
$\Box$

Because the minimum in (\ref{4.1}) is achieved, we can sharpen Lemma 3.3 
as follows.

{\bf Lemma 4.3.} {\em Given any two nodes $x$ and $y$ in ${\cal M}$,
there exists a path $Q_{xy}$ terminating at $x$ and $y$ for which 
$|Q_{xy}|=d(x,y)$.}

There may be more than one such path.  We call each of them an
{\em $x$-to-$y$ geodesic}.

\vspace{.07in}
5. Ordinals and Ranks

As we have seen, the distance between any two nodes of $\cal M$ is
a countable ordinal.  However, given any $x\in{\cal M}$, the set 
$\{d(x,y)\!: y\in{\cal M}\}$ may have no maximum.  For example, 
this is the case for a one-ended 0-path $P^{0}$ where $x$ is any 
fixed node of $P^{0}$ and $y$ ranges through all the 0-nodes of 
$P^{0}$.  On the other hand, for finite graphs the said maximum 
exists and is the ``eccentricity'' of $x$.  We will be able to define an 
``eccentricity'' for every node of $\cal M$ if we expand the set 
$\aleph_{1}$ of countable ordinals into the set $\cal R$ of ranks
\cite[page 4]{tgen}, \cite[page 4]{pg}.  This is done by inserting an 
{\em arrow rank} $\vec{\rho}$ immediately before each limit-ordinal rank
$\rho\in\aleph_{1}$.  $\cal R$ looks like\footnote{As was done in the prior 
works, we treat 0 as the first limit ordinal and $\vec{0}$ as the 
first arrow rank, but in this paper $\vec{0}$ will never be used.  
All our arrow ranks will be understood to be other than $\vec{0}$.}
\[ {\cal R}\;=\;\{\vec{0},0,1,2,\ldots,\vec{\omega},\omega,\omega+1,
\ldots,\vec{\omega\cdot 2},\omega\cdot 2,\omega\cdot 2 +1,\ldots,
\vec{\omega\cdot n},\omega\cdot n, \omega\cdot n+1,\ldots, \]
\[ \vec{\omega^{2}},\omega^{2},\omega^{2}+1,\ldots,\vec{\omega^{k}},\omega^{k},
\omega^{k}+1, \ldots,\vec{\omega^{\omega}},\omega^{\omega},\omega^{\omega}+1,\ldots\}. \]
Note that the set of all ranks is well-ordered. Indeed, there is an 
order-preserving bijection from $\cal R$ to $\aleph_{1}$ obtained 
by replacing each rank by its successor rank.  Since $\aleph_{1}$
is well-ordered, so, too, is $\cal R$.

In accordance with two Aristotelian ideas \cite[page 3]{ru},
we can view each transfinite (successor or limit) ordinal as an 
``actual infinity'' because distances between nodes can assume
those values,
whereas each arrow rank (other than $\vec{0}$) can be viewed as a 
``potential infinity'' because distances can only increase toward
and approach an arrow rank without achieving it.

The arrow ranks served as a notational convenience in the prior works
\cite{tgen} and \cite {pg}, but, for the purposes of this paper, 
we wish to define arrow ranks 
in terms of sequences of countable ordinals.

Let ${\cal A}$ be any set of countable ordinals having a countable 
ordinal $\zeta$
as an upper bound (i.e., $\zeta\geq\alpha$ for all $\alpha\in {\cal A}$).
Let ${\cal D}$ be the set of countable ordinals, each of which is greater than every 
member of ${\cal A}$ and is no greater than $\zeta$.  If ${\cal D}$ is empty, 
${\cal A}$ has a greatest member, namely, $\zeta$. 
So, assume ${\cal D}$ is not empty.  By well-ordering, ${\cal D}$
has a least member $\lambda$.  If $\lambda$ is a successor ordinal,
then ${\cal A}$ has a greatest member, namely, $\lambda-1$;
in this case, $\lambda-1$ is either a successor ordinal or a limit ordinal.
We also denote 
$\lambda -1$ by $\sup {\cal A}$.
If $\lambda$ is a limit ordinal, then there exists an increasing sequence
$\{ \alpha_{k}\}_{k=0}^{\infty}$ contained in ${\cal A}$
such that, for each $\gamma\in {\cal A}$, $\alpha_{k}>\gamma$
for all $k$ sufficiently large (i.e., there exists a $k_{0}$ such 
that $\alpha_{k}>\gamma$ for all $k\geq k_{0}$).  

With $\lambda$ being a nonzero limit ordinal, we define the arrow-rank 
$\vec{\lambda}$ as an equivalence class of such increasing sequences
of ordinals, where two such sequences $\{\alpha_{k}\}_{k=0}^{\infty}$
and $\{\beta_{k}\}_{k=0}^{\infty}$ (not necessarily in $\cal A$ now)
are taken to be equivalent if, for each 
$\gamma$ less than $\lambda$, there exists a natural number $k_{0}$
such that $\gamma<\alpha_{k},\beta_{k}<\lambda$ for all $k>k_{0}$.
The axioms of an equivalence relationship are clearly satisfied.
Each such sequence $\{\alpha_{k}\}_{k=0}^{\infty}$ is a 
representative of $\vec{\lambda}$, and we say that 
$\{ \alpha_{k}\}_{k=0}^{\infty}$ {\em reaches} $\vec{\lambda}$.  
In this case, we let $\sup{\cal A}$ denote $\vec{\lambda}$.

Note that this equivalence class of increasing sequences is different 
from the set of ordinals less than $\lambda$.  The latter is $\lambda$
itself by the definition of ordinals.
Also, we are distinguishing this equivalence class from the limit $\lambda$
of any such sequence in the equivalence class 
\cite[pages 165-166]{ab}.

We summarize these definitions as follows.

{\bf Lemma 5.1.}  {\em If ${\cal A}$ is any set of 
countable ordinals that are bounded above by a countable ordinal
$\zeta$, then $\sup{\cal A}$ exists either as a (successor or limit) ordinal or as an 
arrow rank.}

\vspace{.07in}
6. Eccentricities and Related Ideas

First of all, note that the lengths of all paths in a $\nu$-graph
${\cal G}^{\nu}$ are bounded by $\omega^{\nu+1}\cdot 2$ because the longest possible paths
in ${\cal G}^{\nu}$ are the endless paths of rank $\nu$.  Therefore,
all distances in ${\cal G}^{\nu}$ are also bounded above by $\omega^{\nu+1}\cdot 2$.

The eccentricity $e(x)$ of any node $x\in {\cal M}$ is 
defined by
\begin{equation}
e(x)\;=\; \sup\{ d(x,y)\!: y\in {\cal M}\}.  \label{6.1}
\end{equation}
Two cases arise:  First, the supremum is achieved at some node
$\hat{y}\in {\cal M}$.  In this case, $e(x)$ is an ordinal; so, we can 
replace ``$\sup$'' by ``$\max$'' in (\ref{6.1}) and write
$e(x)=d(x,\hat{y})$.  Second, the supremum is not achieved at
any node in ${\cal M}$.  In this case, $e(x)$ is an arrow rank.

The ideas of radii and diameters for finite graphs \cite[page 32]{b-h},
\cite[page 21]{c-l} can also be extended transfinitely.
Given ${\cal G}^{\nu}$ and ${\cal M}$, the {\em radius} rad$({\cal G}^{\nu},{\cal M})$ 
is the least eccentricity among the nodes of ${\cal M}$:
\begin{equation}
{\rm rad}({\cal G}^{\nu},{\cal M})\;=\;\min\{ e(x)\!: x\in {\cal M}\}  \label{6.2}
\end{equation}
We also denote this simply by rad with the understanding that
${\cal G}^{\nu}$ and ${\cal M}$ are given.  The minimum exists 
as a rank (either as on ordinal 
or as an arrow rank) because the set of ranks is well-ordered.
Thus, there will be at least one $x\in {\cal M}$ with 
$e(x)=\,$rad.  

Furthermore, the {\em diameter} diam$({\cal G}^{\nu},{\cal M})$ 
is defined by 
\begin{equation}
{\rm diam}({\cal G}^{\nu},{\cal M})\;=\; 
\sup_{x\in {\cal M}} \sup_{y\in {\cal M}} \{d(x,y)\}\;=\;\sup\{ d(x,y)\!:
x,y\in {\cal M}\}. \label{6.3}
\end{equation}
With ${\cal G}^{\nu}$ and ${\cal M}$ understood, we denote the 
diameter simply by diam. 
As we have noted before, each $d(x,y)$ is no greater than 
$\omega^{\nu+1}\cdot 2$.  So, by Lemma 5.1, 
diam exists
either as an ordinal or as an arrow rank.  In effect, 
we are defining diam as 
the ``largest'' of the eccentricities.

The ideas of the center and periphery of finite graphs can also be extended.
The {\em center} of $({\cal G}^{\nu},{\cal M})$ 
is the set of nodes in ${\cal M}$
having the least eccentricity, namely, rad. The center is never empty.  

The {\em periphery}
of $({\cal G}^{\nu},{\cal M})$ is the 
set of nodes in ${\cal M}$ having the greatest 
eccentricity, namely, diam.  If diam is an ordinal, there will be
at least two nodes of ${\cal M}$ in the periphery.  Indeed, if there
did not exist at least two nodes in the periphery, then the ordinal
diam could only be approached from below by distances between
pairs of nodes that are less than diam;  thus, the supremum 
in (\ref{6.3}) would have to be an arrow rank---a contradiction.
On the other hand, if diam is an arrow rank, the
periphery can have any positive number of nodes---even just one or 
an infinity of them, as the 
following examples will show.  It seems that the periphery will never be 
empty, but presently this is only a conjecture. 

{\bf Example 6.1.} Let ${\cal G}^{0}$ be a one-ended 0-path 
with ${\cal M}$ being the set of all 0-nodes.  (We do not assign a 
1-node at the path's infinite extremity.)  Then, every 0-node has an 
eccentricity of $\vec{\omega}$.  Thus, rad = diam = $\vec{\omega}$, and 
${\cal M}$ is both the center and the periphery of 
$({\cal G}^{0},{\cal M})$. 

We have here a graphical realization of Aristotle's potential infinity
represented by the eccentricity $\vec{\omega}$.  However, this 
potential infinity can be made into an actual infinity by appending
a 1-node at the infinite extremity of this 0-path.  In Examples 
6.4 and 6.5 below, it will not be possible to convert the potential 
infinities therein into actual infinities.
$\Box$

{\bf Example 6.2.} Consider the 1-graph of Fig. 1 with ${\cal M}$ 
being the set of all 0-nodes along with $x_{a}^{1}$ and $x_{b}^{1}$.
(Ignore $x_{abc}^{1}$ and all other 1-nodes.)  The 0-nodes to 
the left of the 1-nodes all have the eccentricity $\omega+1$.
Also, $e(x_{a}^{1})=e(x_{b}^{1})=\omega\cdot 2 +1$, and
$e(y_{a}^{0})=e(y_{b}^{0})=\omega\cdot 2 +2$.
Thus, rad = $\omega+1$ and diam = $\omega\cdot 2 +2$.  The center
consists of all the 0-nodes to the left of the 1-nodes, and the 
periphery is $\{ y_{a}^{0},y_{b}^{0}\}$. By appending more ``end'' branches 
incident to either $x_{a}^{1}$ or $x_{b}^{1}$, we can increase 
the number of nodes in the periphery.
$\Box$

{\bf Example 6.3.}  Now, consider the 1-graph obtained from
Fig. 1 by deleting the branches $d$ and $e$ and the 0-nodes 
$y_{a}^{0}$ and $y_{b}^{0}$ but appending a new branch incident
to $x_{1}^{0}$ and $x_{a}^{1}$.  Let ${\cal M}$ be all the nodes.
Then, the eccentricity of every node is $\omega$.  Thus, 
rad = diam = $\omega$, and the center and periphery are the same, 
namely, ${\cal M}$. $\Box$

{\bf Example 6.4.}  This time, let ${\cal G}^{\nu}$ consist of a one-ended 0-path 
$P_{a}^{0}$ and an endless 0-path $P_{b}^{0}$ forming a 1-loop, 
as shown in Fig. 2.  ${\cal M}$ is now the set of all (maximal) nodes.
$P_{a}^{0}$ starts at the nonmaximal 0-node $z^{0}$ embraced 
by the 1-node $x^{1}$ and reaches the 1-node $y^{1}$.
$P_{b}^{0}$ reaches both $x^{1}$ and $y^{1}$.  Let
$v^{0}$ be a 0-node of $P_{a}^{0}$ at a distance of $k$ from $x^{1}$.
(For $z^{0}$, $k=0$.)  The shortest distance between $v^{0}$
and any node $w^{0}$ of $P_{b}^{0}$ is provided by a path 
that passes through $x^{1}$; it has the length $\omega+k$.
(The path passing through $y^{1}$ and terminating at $v^{0}$ 
and $\omega^{0}$ has length $\omega\cdot 2$.)  
Thus, $e(v^{0})=\omega+k$.  On the other hand,
$e(w^{0})=\vec{\omega\cdot 2}$;  indeed,
$d(w^{0},v^{0})=\omega+k$, which increases indefinitely 
but never achieves $\omega\cdot 2$ as $v^{0}$ approaches 
$y^{1}$.  Furthermore, $e(x^{1})=e(y^{1})=\omega$.
Thus, rad = $\omega$, diam = $\vec{\omega\cdot 2}$, the center
is $\{x^{1}, y^{1})$, and the periphery is the set
of all the 0-nodes of $P_{b}^{0}$. 

Here is another graphical realization of Aristotle's potential infinity,
this time one represented by the eccentricity $\vec{\omega\cdot 2}$.
In contrast to Example 6.1, there is no way of appending another
transfinite node in order to convert this potential infinity into an 
actual infinity $\omega\cdot 2$. Thus, we have here an incontrovertible 
representation of Aristotle's potential infinity.
$\Box$

{\bf Example 6.5.}  Consider now the 1-graph of Fig. 3.  The eccentricities 
of the nodes are as follows:  $e(x_{k})=\omega\cdot 2+k$ for
$k=1,2,3,\ldots\;$; $e(y^{1})=\omega\cdot 2$; $e(z_{k})^{0}
=\vec{\omega\cdot 2}$ for
$k=\ldots, -1,0-1,\ldots\;$; $e(w^{1})=\vec{\omega\cdot 3}$.
Thus, rad $=\vec{\omega\cdot 2}$, diam $=\vec{\omega\cdot 3}$, 
the center is $\{z_{k}\!:k=\ldots,-1,0,1,\ldots\}$,
and the periphery is the singleton $\{w^{1}\}$. 

Here, we have two different representations of Aristotle's potential infinity,
given by the eccentricities $\vec{\omega\cdot 2}$ and $\vec{\omega\cdot 3}$;
these, too, cannot be converted into ordinals by appending 
transfinite nodes.  Note also that the 
periphery has only one node---in contrast to the peripheries of finite graphs, 
which must have two or more nodes.
$\Box$ 

These examples can immediately be converted into examples for 
graphs of higher ranks by replacing branches by endless paths, all
of the same rank.  For instance, if every branch is replaced
by an endless path of rank $\nu-2$, then every 0-node becomes 
a $(\nu-1)$-node, and every 1-node becomes a $\nu$-node.
Of course, there are far more
complicated $\nu$-graphs.

\vspace{.07in}
7. Some General Results

Henceforth, let it be understood that the metrizable set $\cal M$ 
of nodes has been chosen and fixed for the $\nu$-graph ${\cal G}^{\nu}$
at hand and that any node we refer to is in $\cal M$.

For any rank $\rho$ with $0\leq \rho\leq\nu$, a $\rho$-section
${\cal S}^{\rho}$ of ${\cal G}^{\nu}$ is defined as the subgraph
of ${\cal G}^{\nu}$ induced by a maximal set of branches that are 
$\rho$-connected.\footnote{This is the 
same definition of a $\rho$-section as that
given in \cite[page 49]{tgen} but is somewhat more general than that
of \cite[page 36]{pg}.}  By virtue of Condition 3.2, the $\rho$-sections 
partition ${\cal G}^{\nu}$ (that is, each branch is in one and only 
one $\rho$-section) \cite[Corollary 3.5-6]{tgen}.

We now define a bordering node of ${\cal S}^{\rho}$ to be a node 
of rank larger than $\rho$ that is incident to ${\cal S}^{\rho}$.
Thus, the bordering node embraces as $\alpha$-tip 
$(\alpha\leq \rho)$ traversed by ${\cal S}^{\rho}$;  in other words,
there is a one-ended $\alpha$-path in ${\cal S}^{\rho}$ whose
$\alpha$-tip is embraced by the bordering node.  Also, we define
an internal node of ${\cal S}^{\rho}$ to be a (maximal) node of rank $\rho$
or less contained in ${\cal S}^{\rho}$.

The idea of a component is similar to but different from 
a $(\nu-1)$-section.  A component of a subgraph ${\cal H}$
of ${\cal G}^{\nu}$ is a subgraph of ${\cal H}$ induced by a maximal set of 
branches in ${\cal H}$ that are $\nu$-connected \cite[page 49]{tgen}.
Because ${\cal G}^{\nu}$ is $\nu$-connected, it has just one component, 
namely, itself.
However, a proper subgraph ${\cal H}$ of ${\cal G}^{\nu}$ may have
many components.  For example, if ${\cal H}$ consists of two $(\nu-1)$-sections
that do not share any bordering nodes, then each of them is a 
component of ${\cal H}$.

In the next theorem, ${\cal S}^{\rho}$ is any $\rho$-section 
whose bordering nodes are incident to 
${\cal S}^{\rho}$ only through $\rho$-tips.  
In Fig. 1, $x_{a}^{1}$ and $x_{b}^{1}$ are bordering 
nodes of the 0-section to the left of those nodes, and the condition
is satisfied, that is, those 1-nodes
are incident 
to that 0-section only through 0-tips.  However, branch $d$ induces a 
0-section by itself, and the condition is not satisfied because
$d$ reaches $x_{a}^{1}$ through a $(-1)$-tip
(i.e., a tip of branch $d$);  similarly for $e$ and 
$x_{b}^{1}$.  In Fig. 2, $P_{a}^{0}$ and $P_{b}^{0}$ are 
different 0-sections;  $P_{b}^{0}$ satisfies the condition, but
$P_{a}^{0}$ does not because of node $x^{1}$ and its embraced $(-1)$-tip.

{\bf Theorem 7.1.}  {\em Let ${\cal S}^{\rho}$ be a $\rho$-section 
in ${\cal G}^{\nu}$
$(0\leq \rho< \nu)$ all of whose bordering nodes are incident to 
${\cal S}^{\rho}$ only through $\rho$-tips.  
Then, all the internal nodes of ${\cal S}^{\rho}$
have the same eccentricity.}

{\em Proof.}  By virtue of our hypothesis and the $\rho$-connectedness
of ${\cal S}^{\rho}$, for any internal node $x^{\alpha}$ $(\alpha\leq \rho)$
and any bordering node $z^{\gamma}$
$(\gamma >\rho)$ of ${\cal S}^{\rho}$ in ${\cal M}$, there is a 
representative $\rho$-path 
$P^{\rho}$ for a $\rho$-tip embraced by $z^{\gamma}$ and lying in 
${\cal S}^{\rho}$, and 
there also is a two-ended path $Q$ lying in ${\cal S}^{\rho}$
and terminating at $x^{\alpha}$ and a node of $P^{\rho}$.
So, by Condition 3.2 and 
\cite[Corollary 3.5-4]{tgen}, there is in $P\cup Q$ a one-ended
$\rho$-path $R^{\rho}$ that terminates at $x^{\alpha}$ and reaches
$z^{\gamma}$ through a $\rho$-tip.  Moreover, all paths that terminate at 
$x^{\alpha}$, that lie in ${\cal S}^{\rho}$, 
and that reach $z^{\gamma}$ must be one-ended 
$\rho$-paths. Therefore, $d(x^{\alpha},z^{\gamma})=\omega^{\rho+1}$.
For any other node $y^{\beta}$ $(\beta\leq \rho)$ in ${\cal S}^{\rho}$, we 
have $d(x^{\alpha},y^{\beta})< \omega^{\rho+1}$ by the 
$\rho$-connectedness of 
${\cal S}^{\rho}$.  So, if ${\cal G}^{\nu}$ consists only 
of ${\cal S}^{\rho}$ and its bordering nodes (so that $\nu=\rho+1$),
we can conclude that 
$e(x^{\alpha})=\omega^{\rho+1}$, whatever be the choice of the internal node 
$x^{\alpha}$ in ${\cal S}^{\rho}$ and in ${\cal M}$.

Next, assume that there is a node 
$v^{\delta}$ of ${\cal G}^{\nu}$ 
lying outside of ${\cal S}^{\rho}$ 
and different from all the bordering nodes of ${\cal S}^{\rho}$.
By the $\nu$-connectedness of ${\cal G}^{\nu}$,
there is a path $P_{xv}$ terminating at $x^{\alpha}$ and $v^{\delta}$.
Let $z^{\gamma}$ now be the last bordering node of ${\cal S}^{\rho}$
that $P_{xv}$ meets.  Let $P_{zv}$ be that part of $P_{xv}$
lying outside of ${\cal S}^{\rho}$.  Then, 
by what we have shown above, there is a 
one-ended $\rho$-path $Q^{\rho}_{xz}$ that terminates 
at $x^{\alpha}$, lies in ${\cal S}^{\rho}$, and
reaches $z^{\gamma}$ through a $\rho$-tip.  Then, $R_{xv}=Q_{xz}\cup P_{zv}$
is a two-ended path that terminates at $x^{\alpha}$ and $v^{\delta}$.
Moreover, $|R_{xv}|\leq|P_{xv}|$.

Now, let $y^{\beta}$ $(\beta\leq \rho)$ be any other internal node of 
${\cal S}^{\rho}$ in ${\cal M}$ (i.e., different from $x^{\alpha}$).
Again, there is a one-ended $\rho$-path $Q_{yz}^{\rho}$ satisfying
the same conditions as $Q_{xz}^{\rho}$.  We have 
$d(x^{\alpha},z^{\gamma})=d(y^{\beta},z^{\gamma})=\omega^{\rho+1}$.
Let $R_{yv}=Q_{yz}\cup P_{zv}$.  Thus, $|R_{xv}|=|R_{yv}|$.
We have shown that, for each one-ended path $R_{xv}$
terminating at $x^{\alpha}$ and $v^{\delta}$ and passing through exactly
one bordering node of $z^{\gamma}$ of ${\cal S}^{\rho}$,
there is another path $R_{yv}$ of the same length terminating at 
$y^{\beta}$ and $v^{\delta}$ and identical to $R_{xv}$ outside 
${\cal S}^{\rho}$.  
It follows that $d(x^{\alpha},v^{\delta})=d(y^{\beta},v^{\delta})$.
We can conclude that $e(x^{\alpha})=e(y^{\beta})$ whatever be the choices of 
$x^{\alpha}$ and $y^{\beta}$ in ${\cal S}^{\rho}$ and ${\cal M}$.  $\Box$

Figs. 1, 2, and 3 provide examples for Theorem 7.1.  In Fig. 1,
all the 0-nodes to the left of the 1-nodes have the same eccentricity
$\omega+1$ in accordance with the theorem.
In Fig. 2, all the nodes of $P_{b}^{0}$ have the same eccentricity
$\vec{\omega\cdot 2}$, whereas the eccentricities of the nodes of 
$P_{a}^{0}$ vary;  this, too, conforms with Theorem 7.1.
Similarly, in Fig. 3, the nodes $z_{k}^{0}$ have the same eccentricities,
but the nodes $x_{k}^{0}$ have differing eccentricities.

A standard result \cite[page 21]{c-l} can be extended to the 
transfinite case, albeit in a more complicated way.  Given
${\cal G}^{\nu}$ and $\cal M$, rad may be either an ordinal or an arrow rank.  
If it is an arrow rank, we let rad$^{+}$ denote the limit ordinal immediately 
following rad.

{\bf Theorem 7.2.}  {\em  
\begin{description}
\item[(i)] If rad is an ordinal, then rad $\leq$ diam $\leq$ rad$\cdot 2$.
\item[(ii)] If rad is an arrow rank, then rad $\leq$ diam $\leq$
rad$^{+}\cdot 2$.
\end{description}
}

{\bf Proof.}  The proofs of (i) and (ii) are much the same.  So, let us
consider (ii) alone.  That rad $\leq$ diam follows directly from the 
definitions (\ref{6.1}), (\ref{6.2}), and (\ref{6.3}).
Next, by the definition 
of the diameter (\ref{6.3}), we can choose two sequences
$\{y_{k}\}_{k=0}^{\infty}$ and $\{z_{k}\}_{k=0}^{\infty}$
of nodes 
such that the sequence $\{ d(y_{k},z_{k})\}_{k=0}^{\infty}$
approaches or achieves diam.  Let $x$ be any node in the center.
By the triangle inequality,
\[ d(y_{k},z_{k})\;\leq\; d(y_{k},x)\,+\, d(x,z_{k}). \]
Now, $d(y_{k},x)\leq\;$rad $\leq$ rad$^{+}$, and similarly for 
$d(x,z_{k})$.  Therefore, 
\[d(y_{k},z_{k})\leq\; {\rm rad}^{+}\;+\;{\rm rad}^{+} \;=
{\rm rad}^{+}\;\cdot 2. \]  
$\Box$

Another standard result is that the nodes of any finite graph comprise 
the center of some finite connected graph \cite[page 22]{c-l}.\footnote{This
result extends immediately to infinite 0-graphs with infinitely many 0-nodes.
We are now considering transfinite graphs of ranks 1 or greater.}  
This, too, 
can be extended transfinitely---in fact, in several ways, but 
the proofs are more complicated than that for finite graphs.
Nonetheless, the scheme of the proofs remains the same.
First, we need the following lemma.  In the following, 
$\nu-1$ denotes $\vec\omega$
when $\nu=\omega$. 

{\bf Lemma 7.3.}  {\em  Let 
${\cal S}^{\nu-1}$ be a $(\nu-1)$-section of ${\cal G}^{\nu}$, where 
$1\leq \nu\leq\omega$ and $\nu\neq\vec{\omega}$.
Let $u^{\nu}$ 
be a $\nu$-node incident to ${\cal S}^{\nu-1}$ (thus, a 
bordering node of ${\cal S}^{\nu-1}$), and
let $x^{\alpha}$ $(\alpha<\nu)$ be 
an $\alpha$-node in ${\cal S}^{\nu-1}$
(thus, an internal node of ${\cal S}^{\nu-1}$).
Then, there exists in ${\cal S}^{\nu-1}$
a two-ended path of length no larger than 
$\omega^{\nu}$ connecting $x^{\alpha}$ and $u^{\nu}$.}

{\bf Proof.}  That $u^{\nu}$ is incident to ${\cal S}^{\nu-1}$
means that there is in ${\cal S}^{\nu-1}$ a one-ended $\beta$-path 
$P^{\beta}$ with 
$\beta\leq \nu-1$ whose $\beta$-tip is embraced 
by $u^{\nu}$.  Let $P^{\beta+1}$ be the two-ended path obtained
by appending to $P^{\beta}$ the $(\beta+1)$-node $y^{\beta+1}$
embraced by $u^{\nu}$ and reached by $P^{\beta}$.  
($y^{\beta+1}$ will not be maximal if $\beta+1<\nu$; 
otherwise, $y^{\beta+1}=u^{\nu}$.)  The length $|P^{\beta+1}|$
of $P^{\beta+1}$ is equal to 
$\omega^{\beta+1}$ because $P^{\beta+1}$ 
traverses only one $\beta$-tip;  all other tips traversed by 
$P^{\beta+1}$ are of lesser rank.  Let $z^{\gamma}$ be any node
(not necessarily maximal) of 
$P^{\beta}$; thus, $\gamma\leq \beta$.  By the $(\nu-1)$-connectedness
of ${\cal S}^{\nu-1}$, 
there is in ${\cal S}^{\nu-1}$ a two-ended $\lambda$-path $Q^{\lambda}$
$(0\leq \lambda\leq \nu-1)$ terminating at $x^{\alpha}$ and $z^{\gamma}$.
The tips traversed by $Q^{\lambda}$ have ranks no greater than 
$\lambda-1$, hence, no greater 
than $\nu-2$.  By \cite[Corollary 3.5-4]{tgen}, 
there is a two-ended path $R^{\delta}$ in $P^{\beta+1}\cup Q^{\lambda}$
terminating at $x^{\alpha}$ and $y^{\beta+1}$.  All the tips traversed by 
$R^{\delta}$ are of ranks no greater than $\nu-1$, 
and there is at most one traversed tip of rank $\nu-1$.
Hence, the length of $R^{\delta}$ satisfies
$|R^{\delta}|\leq \omega^{\nu}$.  $\Box$

Given any $\nu$-graph ${\cal G}^{\nu}$ with 
$1\leq \nu\leq\omega$ and $\nu\neq\vec{\omega}$, 
let us construct a larger $\nu$-graph
${\cal H}^{\nu}$ by appending six additional $\nu$-nodes $p_{i}^{\nu}$ and 
$q_{i}^{\nu}$ $(i=1,2,3)$ and also appending isolated endless 
$(\nu-1)$-paths\footnote{An isolated endless path embraces 
no tips other than the ones it 
traverses.  Thus, to reach any other part of a graph in which
the isolated path is a subgraph, one must proceed through a terminal tip
of that path.} that reach $\nu$-nodes as shown in Fig. 4.  
Such paths connect $p_{1}^{\nu}$ to $p_{2}^{\nu}$, $p_{2}^{\nu}$ 
to $p_{3}^{\nu}$, $p_{3}^{\nu}$ to every $\nu$-node in ${\cal G}^{\nu}$, and 
similarly for $p_{i}^{\nu}$ replaced by $q_{i}^{\nu}$.
Note that the singleton end-nodes $p_{1}^{\nu}$ and $q_{1}^{\nu}$
can be included in the chosen metrizable set ${\cal M}$ for ${\cal H}^{\nu}$.  
All the other $\nu$-nodes of ${\cal H}^{\nu}$ are nonsingletons
and therefore are in ${\cal M}$, too.

{\bf Theorem 7.4.}  {\em The $\nu$-nodes 
of ${\cal G}^{\nu}$ ($1\leq \nu\leq\omega$, $\nu\neq\vec{\omega}$) 
comprise the center of 
${\cal H}^{\nu}$, and the periphery of ${\cal H}^{\nu}$ 
is $\{p_{1}^{\nu},q_{1}^{\nu}\}$.}

{\bf Proof.}  We look for bounds on the eccentricities of all the nodes 
in ${\cal M}$.  Let $x^{\alpha}$ and 
$y^{\beta}$ be any two nodes  
whose ranks satisfy 
$0\leq \alpha,\beta <\nu$. It follows that $x^{\alpha}$ 
(resp. $y^{\beta}$) is an internal node of a $(\nu-1)$-section in 
${\cal G}^{\nu}$, and that section has at least one $\nu$-node $u^{\nu}$
(resp. $v^{\nu}$) as a bordering node because ${\cal G}^{\nu}$
is $\nu$-connected.  By the triangle inequality,
\[ d(x^{\alpha},y^{\beta})\;\leq\; d(x^{\alpha},u^{\nu})\,+\,
d(u^{\nu},p_{3}^{\nu})\,+\,d(p_{3}^{\nu},v^{\nu})\,+\,
d(v^{\nu},y^{\beta}). \]
By Lemma 7.3, $d(x^{\alpha},u^{\nu})\leq \omega^{\nu}$ and 
$d(v^{\nu},y^{\beta})\leq\omega^{\nu}$.  
Clearly, $d(u^{\nu},p_{3}^{\nu})=
d(p_{3}^{\nu},v^{\nu})=\omega^{\nu}\cdot 2$.  Thus, 
$d(x^{\alpha},y^{\beta})\leq \omega^{\nu}\cdot 6$.  This also shows
that, for any $\nu$-node $v^{\nu}$ in ${\cal G}^{\nu}$,
$d(x^{\alpha},v^{\nu})\leq \omega^{\nu}\cdot 5$.
Since $d(u^{\nu},p_{1}^{\nu})=\omega^{\nu}\cdot 6$, 
we have $d(x^{\alpha},p_{1}^{\nu})=
d(x^{\alpha},u^{\nu})+d(u^{\nu},p_{1}^{\nu})
\leq \omega^{\nu}+\omega^{\nu}\cdot 6=\omega^{\nu}\cdot 7$.
Now, $d(x^{\alpha},u^{\nu})\geq 1$ because there is at least one branch in
any path connecting $x^{\alpha}$ and $u^{\nu}$.  Thus, 
we also have $d(x^{\alpha},p_{1}^{\nu})\geq \omega^{\nu}\cdot 6+1$.  
Note also that the distance from $x^{\alpha}$ to any node of the 
appended endless paths is strictly less than $\omega^{\nu}\cdot 7$.
All these results hold for $p_{i}^{\nu}$ replaced by $q_{i}^{\nu}$.
Altogether then, we can conclude the following:  For any 
node in ${\cal G}^{\nu}$ of 
rank less than $\nu$, say, $x^{\alpha}$, the eccentricity $e(x^{\alpha})$ 
of $x^{\alpha}$ is bounded as follows:
\[ \omega^{\nu}\cdot 6+1\;\leq\; e(x^{\alpha})\;\leq\; 
\omega^{\nu}\cdot 7 \]

Next, consider any two $\nu$-nodes of ${\cal G}^{\nu}$, say, $u^{\nu}$
and $v^{\nu}$ again.  By what we have already shown, 
$d(u^{\nu},v^{\nu})\leq\omega^{\nu}\cdot 4$, and 
$d(u^{\nu},p_{1}^{\nu})=d(u^{\nu},q_{1}^{\nu})=\omega^{\nu}\cdot 6$.
The distance from $u^{\nu}$ to any node of the appended endless
$(\nu-1)$-paths is less than $\omega^{\nu}\cdot 6$.
Also, for any node $y^{\beta}$ in ${\cal G}^{\nu}$ of rank less than $\nu$,
$d(u^{\nu},y^{\beta})\leq \omega^{\nu}\cdot 5$.  So, 
the largest distance between $u^{\nu}$ and any other node in
${\cal H}^{\nu}$ is equal to $\omega^{\nu}\cdot 6$;  that is, 
$e(u^{\nu})=\omega^{\nu}\cdot 6$.

Finally, we have $e(p_{3}^{\nu})=e(q_{3}^{\nu})=\omega^{\nu}\cdot 8$,
$e(p_{2}^{\nu})=e(q_{2}^{\nu})=\omega^{\nu}\cdot 10$,
and $e(p_{1}^{\nu})=e(q_{1}^{\nu})=\omega^{\nu}\cdot 12$. 
The eccentricities of the nodes of the appended endless paths lie between 
these values.

We have considered all cases.  Comparing these equalities and inequalities
for all the eccentricities, we can draw the conclusion of the theorem.
$\Box$

As an immediate corollary, we have the following generalization 
of a result for finite graphs.

{\bf Corollary 7.5.}  {\em The $\nu$-nodes of ${\cal G}^{\nu}$ 
($1\leq \nu\leq\omega$, $\nu\neq\vec{\omega}$) comprise the 
center of some $\nu$-graph ${\cal H}^{\nu}$.}

Variations of Corollary 7.5 can also be established through 
much the same proofs.  For instance, all the nodes of 
${\cal G}^{\nu}$ of one or more specified ranks can be 
made to comprise the center of some $\nu$-graph.  This is
because the $(\nu-1)$-sections of ${\cal G}^{\nu}$ partition ${\cal G}^{\nu}$.
Still more generally, if ${\cal G}^{\nu}$ has only finitely many $\nu$-nodes,
any arbitrary set of nodes of ${\cal G}^{\nu}$ 
in $\cal M$ can be made the center
simply by appending enough endless $(\nu-1)$-paths in series.

\vspace{.07in}
8. When All the Nodes of Highest Rank Are Pristine

The nodes of highest rank in ${\cal G}^{\nu}$ are the $\nu$-nodes, of course.
A $\nu$-node is said to be pristine if it does not embrace a node of 
lower rank.  Thus, a pristine $\nu$-node consists only of $(\nu-1)$-tips.
Henceforth, we assume the following. 

{\bf Condition 8.1.}  {\em All the $\nu$-nodes are pristine.}

Because of this, we can view the $\nu$-nodes as lying only at infinite
extremities of the $(\nu-1)$-sections to which they are incident
because they can be reached only through $(\nu-1)$-tips
of such sections; they are, in fact, 
the bordering nodes of such $(\nu-1)$-sections.
All the other nodes are of ranks less than $\nu$ and are internal nodes of 
$(\nu-1)$-sections. 
Given any $(\nu-1)$-section ${\cal S}^{\nu-1}$, the set 
of all internal nodes of 
${\cal S}^{\nu-1}$ will be denoted by $i({\cal S}^{\nu-1})$
and will be called the {\em interior} of ${\cal S}^{\nu-1}$.

A boundary $\nu$-node of ${\cal S}^{\nu-1}$ 
is defined to be a bordering $\nu$-node 
that contains $(\nu-1)$-tips of ${\cal S}^{\nu-1}$ and also $(\nu-1)$-tips
of one or more other $(\nu-1)$-sections of ${\cal G}^{\nu}$.
Thus, a boundary node lies at the infinite extremities of two or more 
$(\nu-1)$-sections and thereby connects them.
A $\nu$-path can pass from the interior of one $(\nu-1)$-section into the 
interior of another $(\nu-1)$-section only by passing through a 
boundary node.

Another assumption we henceforth impose is the following.

{\bf Condition 8.2.}  {\em There are only finitely many boundary $\nu$-nodes
throughout ${\cal G}^{\nu}$.} 

Nevertheless, each
$(\nu-1)$-section may have infinitely many incident non-boundary bordering 
$\nu$-nodes, and each boundary $\nu$-node may be incident to infinitely many 
$(\nu-1)$-sections.

{\bf Lemma 8.3.}  {\em  If $P$ is a two-ended $\nu$-path, then
$|P|=\omega^{\nu}\cdot k$, where $k$ is the number of $(\nu-1)$-tips
traversed by $P$.}

{\bf Proof.}  This follows directly from the definition of the 
length $|P|$ and the fact that all $\nu$-nodes are pristine.  $\Box$

Because all $\nu$-nodes are now pristine, we can strengthen 
Lemma 7.3 as follows.

{\bf Lemma 8.4.}  {\em 
\begin{description}
\item[(a)] Let $x^{\nu}$ be a bordering node of a $(\nu-1)$-section 
${\cal S}^{\nu-1}$, and let $z$ be an internal node 
of ${\cal S}^{\nu-1}$.  Then,
there exists a two-ended $\nu$-path $P_{z,x}^{\nu}$
within ${\cal S}^{\nu-1}$ that terminates at $z$ and 
reaches $x^{\nu}$ through its 
one and only $(\nu-1)$-tip.  Moreover, the length $|P_{z,x}^{\nu}|$
of $P_{z,x}^{\nu}$ is $\omega^{\nu}$.  $P_{z,x}^{\nu}$ is a $z$-to-$x$
geodesic.  Finally, all two-ended $\nu$-paths within ${\cal S}^{\nu-1}$
terminating at an internal node of ${\cal S}^{\nu-1}$ and at
a bordering node of ${\cal S}^{\nu-1}$ have the length $\omega^{\nu}$.
\item[(b)] Let $x^{\nu}$ and $y^{\nu}$ be two bordering nodes of 
${\cal S}^{\nu-1}$.  Then, there exists an two-ended 
$\nu$-path $P_{x,y}^{\nu}$
within ${\cal S}^{\nu-1}$ that reaches $x^{\nu}$ and $y^{\nu}$ through
its two $(\nu-1)$-tips.  Moreover, the length $|P_{x,y}^{\nu}|$ 
of $P_{x,y}^{\nu}$ is $\omega^{\nu}\cdot 2$.  $P_{x,y}^{\nu}$ is an
$x^{\nu}$-to-$y^{\nu}$ geodesic.  Finally, all 
two-ended $\nu$-paths in ${\cal S}^{\nu-1}$ reaching two bordering nodes of 
${\cal S}^{\nu-1}$ have the length $\omega^{\nu}\cdot 2$.
\end{description}
}

{\bf Proof.}  The proof of part (a) is much the same as that of 
Lemma 7.3 except that now $P_{z,x}^{\nu}$ is incident to $x^{\nu}$
only through its one and only $(\nu-1)$-tip.

Part (b) is proven similarly, but now we use two representative
one-ended paths, one for each of $x^{\nu}$ and $y^{\nu}$.  $\Box$

{\bf Lemma 8.5.}  {\em Let $x^{\nu}$ be a bordering node of a 
$(\nu-1)$-section ${\cal S}^{\nu-1}$, let $z$ be an internal node of 
${\cal S}^{\nu-1}$, and let $y$ be any node of ${\cal G}^{\nu}$.  
Then, $|d(x^{\nu},y)-d(z,y)|\leq \omega^{\nu}$.}

{\bf Proof.}  By Lemma 8.4(a), $d(z,x^{\nu})=\omega^{\nu}$.  
Since, $d$ is a metric, 
\[ d(x^{\nu},y)\;\leq\;d(x^{\nu},z)\,+\,d(z,y)\;=\;\omega^{\nu}\,+\,d(z,y). \]
Also,
\[ d(z,y)\;\leq\;d(z,x^{\nu})\,+\,d(x^{\nu},y)\;=\;
\omega^{\nu}\,+\,d(x^{\nu},y). \]
These inequalities yield the conclusion.  $\Box$

We will show below (Theorem 8.7) that, as a consequence of Conditions
8.1 and 8.2, no node of ${\cal G}^{\nu}$ can have an arrow-rank eccentricity
and that all eccentricities comprise a finite set of ordinal values.
But, first note that three examples of the occurrence of 
arrow-rank eccentricities are given in Examples 
6.1, 6.4, and 6.5. 
Each of these examples violates either 
Condition 8.1 or Condition 8.2.  

{\bf Theorem 8.6.}  {\em The eccentricities of all the nodes 
are contained within the following finite set of ordinals:
\begin{equation}
\{ \omega^{\nu}\cdot p\!: 1\leq p\leq 2m+2\}  \label{8.3}
\end{equation}
Here, $p$ and $m$ are natural numbers, and $m$ is the number of 
boundary $\nu$-nodes.}

{\bf Proof.}  The eccentricity of any node is at least as large 
as the distance between  any internal node of a $(\nu-1)$-section
and any bordering $\nu$-node of that $(\nu-1)$-section.  Therefore,
Lemma 8.4(a) implies that the eccentricity of any node of ${\cal G}^{\nu}$
is at least $\omega^{\nu}$, whence the lower bound in (\ref{8.3}).
The proof of the upper bound requires more effort.

First of all, we can settle two simple cases by inspection.  
If ${\cal G}^{\nu}$ consists of a single $(\nu-1)$-section with exactly one 
bordering $\nu$-node (in ${\cal M}$, of course), 
then all the nodes of ${\cal G}^{\nu}$
have the eccentricity $\omega^{\nu}$.  If that one and only 
$(\nu-1)$-section for ${\cal G}^{\nu}$ has two or more 
(possibly infinitely many) bordering nodes, the internal nodes have 
eccentricity $\omega^{\nu}$, and the bordering nodes have eccentricity
$\omega^{\nu}\cdot 2$.  In both cases, the conclusion of the theorem
is fulfilled with $m=0$.

We now turn to the general case where ${\cal G}^{\nu}$ has 
at least one boundary 
$\nu$-node and therefore at least two $(\nu-1)$-sections.
${\cal G}^{\nu}$ will have a two-ended $\nu$-path of the following form:
\begin{equation}
P_{0,k}^{\nu}\;=\;\{x_{0},P_{0}^{\nu-1},x_{1}^{\nu},P_{1}^{\nu-1},
\ldots,x_{k-1}^{\nu},P_{k-1}^{\nu-1},x_{k}\}  \label{8.4}
\end{equation}
Because all $\nu$-nodes are pristine, the $x_{i}^{\nu}$ $(i=1,\ldots,k-1)$
are nonsingleton bordering $\nu$-nodes (possibly boundary $\nu$-nodes),
and the $P_{i}^{\nu-1}$ $(i=2,\ldots,k-2)$ are endless $(\nu-1)$-paths.
The same is true of $x_{0}$, $x_{k}$, $P_{1}^{\nu-1}$, and $P_{k-1}^{\nu-1}$
if $x_{0}$ and $x_{k}$ are $\nu$-nodes, too.  If $x_{0}$ 
(resp. $x_{k}$)
is of lower rank, then it is an internal node, and $P_{1}^{\nu-1}$
(resp. $P_{k-1}^{\nu-1}$)
is a one-ended $(\nu-1)$-path.  

Let us first assume that $x_{0}$ and $x_{k}$ are internal nodes
in different $(\nu-1)$-sections.
Let ${\cal S}_{0}^{\nu-1}$ be the $(\nu-1)$-section containing $x_{0}$.  Let 
$x_{i_{1}}^{\nu}$ be the last $\nu$-node in (\ref{8.4})
that is incident to ${\cal S}_{0}^{\nu-1}$.  $x_{i_{1}}^{\nu}$ 
will be a boundary 
$\nu$-node because it is also incident to another $(\nu-1)$-section, say, 
${\cal S}_{1}^{\nu-1}$. If need be,
we can replace the subpath of (\ref{8.4}) between $x_{0}$
and $x_{i_{1}}$ by a $\nu$-path $\{x_{0},Q_{0}^{\nu-1}$,
$x_{i_{1}}^{\nu}\}$, where $Q_{0}^{\nu-1}$ is a one-ended $(\nu-1)$-path 
and resides in ${\cal S}_{0}^{\nu-1}$, to get a shorter overall 
$\nu$-path terminating at 
$x_{0}$ and $x_{k}$.

Now, let ${\cal S}_{1}^{\nu-1}$ be the next $(\nu-1)$-section
after ${\cal S}_{0}^{\nu-1}$ through which our (possibly) reduced path 
proceeds.  Also,
let $x_{i_{2}}^{\nu}$ be the last $\nu$-node in that path
that is incident to ${\cal S}_{1}^{\nu-1}$.  $x_{i_{2}}^{\nu}$ 
will be a boundary 
$\nu$-node incident to ${\cal S}_{1}^{\nu-1}$ and another $(\nu-1)$-section
${\cal S}_{2}^{\nu-1}$.  If need be,
we can replace the subpath between 
$x_{i_{1}}$ and $x_{i_{2}}$ by a $\nu$-path $\{x_{i_{1}}^{\nu}$,
$Q_{1}^{\nu-1},x_{i_{2}}^{\nu}\}$, where $Q_{1}^{\nu-1}$
is an endless $(\nu-1)$-path residing in ${\cal S}_{1}^{\nu-1}$.
This will  yield a still shorter overall $\nu$-path 
terminating at $x_{0}$ and $x_{k}$.

Continuing this way, we will find a boundary $\nu$-node
$x_{i_{j}}^{\nu}$ that is incident to the 
$(\nu-1)$-section ${\cal S}_{j}^{\nu-1}$
containing $x_{k}$.  Finally, we let $Q_{j}^{\nu-1}$ be a one-ended 
$(\nu-1)$-path in ${\cal S}_{j}^{\nu-1}$ terminating 
at $x_{i_{j}}^{\nu}$ and $x_{k}$.  Altogether, we will
have the following two-ended $\nu$-path, which is not longer than 
$P_{0,k}^{\nu}$ (actually shorter if the 
aforementioned replacements were needed).
\begin{equation}
Q_{0,k}^{\nu}\;=\;\{x_{0},Q_{0}^{\nu-1},x_{i_{1}}^{\nu},Q_{1}^{\nu-1},
x_{i_{2}}^{\nu},\ldots,x_{i_{j}}^{\nu},Q_{j}^{\nu-1},x_{k}\}  \label{8.5}
\end{equation}
Because all the $\nu$-nodes herein are boundary nodes and pristine,
the length $|Q_{0,k}|$ is obtained simply by counting the 
$(\nu-1)$-tips traversed by $Q_{0,k}^{\nu}$ and multiplying by 
$\omega^{\nu}$ (see Lemma 8.4).  We get $|Q_{0,k}^{\nu}|=
\omega^{\nu}\cdot(2j)$, where $j\leq m$.
Finally, we note that any geodesic path between $x_{0}$ and $x_{k}$ has 
a length no larger than than $|Q_{0,k}^{\nu}|=\omega^{\nu}\cdot(2j)$.

Next, consider the case where $x_{0}$ is a bordering $\nu$-node of 
${\cal S}_{0}^{\nu-1}$ and $x_{k}$ remains an internal node 
of ${\cal S}_{j}^{\nu-1}$.
$P_{0}^{\nu-1}$ in (\ref{8.4}) will 
be an endless $(\nu-1)$-path residing in some $(\nu-1)$-section
${\cal S}_{0}^{\nu-1}$.  We let $x_{i_{1}}^{\nu}$ be the last $\nu$-node in 
(\ref{8.4}) incident to ${\cal S}_{0}^{\nu-1}$.  Otherwise 
our procedure is as before,
and we can now conclude that $|Q_{0,k}^{\nu}|=\omega^{\nu}\cdot
(2j+1)$ because the passage from $x_{0}$ into ${\cal S}_{0}^{\nu-1}$
traverses one $(\nu-1)$-tip.

The same conclusion, namely, $|Q_{0,k}^{\nu}|=\omega^{\nu}\cdot(2j+1)$
holds if $x_{k}$ is a bordering $\nu$-node and $x_{0}$ is 
an internal node.  Finally, if both $x_{0}$ and $x_{k}$
are bordering $\nu$-nodes, we get $|Q_{0,k}^{\nu}|=
\omega^{\nu}\cdot(2j+2)$.  For all cases, we can assert that the 
geodesic between $x_{0}$ and $x_{k}$ has a length no larger 
than $|Q_{0,k}^{\nu}|=\omega^{\nu}\cdot(2j+2)$.

Now, the eccentricity $e(x_{0})$ for $x_{0}$ is the supremum of the 
lengths of all geodesics starting at $x_{0}$ and terminating at all 
other nodes $x_{k}$.  Since there are only finitely many 
boundary nodes and since every geodesic will have the form 
of (\ref{8.5}), every eccentricity will be a multiple of 
$\omega^{\nu}$ (there are no arrow-rank eccentricities).
Also, since $j\leq m$ where $m$ is the number of 
boundary $\nu$-nodes in ${\cal G}^{\nu}$, we can conclude that $e(x_{0})\leq
\omega^{\nu}\cdot(2m+2)$, whatever be the node $x_{0}$.  $\Box$

The lengths of all
geodesics will reside in the finite set of values (\ref{8.3}).
Consequently, for every node $x_{0}$ of ${\cal G}^{\nu}$ there will be at least 
one geodesic of maximum length starting at $x_{0}$ and terminating 
at some other node $z$ of ${\cal G}^{\nu}$.  Such a geodesic is called an 
{\em eccentric path} for $x_{0}$, and $z$ is called an
{\em eccentric node} for $x_{0}$.  
In general, there are many eccentric paths and eccentric nodes for 
a given $x_{0}$.

{\bf Corollary 8.7.}  {\em Let $x^{\nu}$ be any bordering $\nu$-node of a
$(\nu-1)$-section ${\cal S}^{\nu-1}$ with the eccentricity 
$e(x^{\nu})=\omega^{\nu}\cdot k$, and let $z$ be an internal node of 
${\cal S}^{\nu-1}$ with the eccentricity $e(z)=\omega^{\nu}\cdot p$.
Then, $|k-p|\leq 1$.}

{\bf Proof.}  Let $P_{z,x}$ be the two-ended $\nu$-path 
obtained by appending $x^{\nu}$ to a one-ended $(\nu-1)$-path in 
${\cal S}^{\nu-1}$ that reaches $x^{\nu}$ and terminates at the 
internal node $z$.  By 
Lemma 8.4(a), $P_{z,x}^{\nu}$ is a $z$-to-$x^{\nu}$ geodesic,
and $|P_{z,x}^{\nu}|=d(z,x^{\nu})=\omega^{\nu}$.
Now, let $w$ be any node.  By the triangle inequality for the metric $d$,
\[ d(z,w)\;\leq \; d(z,x^{\nu})\,+\,d(x^{\nu},w)\;=\;
\omega^{\nu}\,+\,d(x^{\nu},w). \]
Next, let $w$ be an eccentric node for $z$.  We get $d(z,w)=e(z)$ and
$e(z)\leq \omega^{\nu}+d(x^{\nu}, w)$.  Moreover, $d(x^{\nu},w)\leq
e(x^{\nu})$.  Therefore, 
\begin{equation}
e(z)\;\leq\;\omega^{\nu}\,+\,e(x^{\nu}).  \label{8.6}
\end{equation}
By a similar argument with $w$ now being an eccentric node for $x^{\nu}$,
we get
\begin{equation}
e(x^{\nu})\;\leq\;\omega^{\nu}\,+\,e(z).  \label{8.7}
\end{equation}
So, with (\ref{8.6}) we have $\omega^{\nu}\cdot p\leq \omega^{\nu}
+\omega^{\nu}\cdot k=\omega^{\nu}\cdot(k+1)$, or $p\leq k+1$.
On the other hand, with (\ref{8.7}) we have in the same way $k\leq p+1$.
Whence our conclusion.  $\Box$

That $k-p$ can equal 0 is verified by the next example.

{\bf Example 8.8.}  Consider the 1-graph of Fig. 5 consisting of a 
one-ended 0-path of 0-nodes $w_{k}^{0}$ and an endless 0-path of 
0-nodes $y_{k}^{0}$ connected in series to two 1-nodes $x^{1}$ and $z^{1}$ 
as shown.  The eccentricities are as follows:
$e(w_{k}^{0})=\omega\cdot 3$ for $k=1,2,3, \ldots\;$,
$e(x^{1})=\omega\cdot 2$, $e(y_{k}^{0})=\omega\cdot 2$ for 
$k=\ldots,-1,0,1,\dots\;$, and $e(z^{1})=\omega\cdot 3$.  Thus,
$e(x^{1})-e(y_{k}^{0})=0$, as asserted.  $\Box$

An immediate consequence of Theorem 8.6 and Corollary 8.7 is the following.

{\bf Corollary 8.9.}  {\em The eccentricities of all the nodes form 
a consecutive set of values in (\ref{8.3}), with the minimum 
(resp. maximum) eccentricity being the radius (resp. diameter) of 
${\cal G}^{\nu}$.}

By virtue of Theorem 7.1 and Condition 8.1, we have that, if 
an ordinal is the eccentricity of an internal node of a 
$(\nu-1)$-section ${\cal S}^{\nu-1}$,
then there will be infinitely many nodes with the same eccentricity,
for example, all the internal nodes of ${\cal S}^{\nu-1}$.  Moreover,
a boundary node of ${\cal S}^{\nu-1}$ may also have that same
eccentricity;  the 1-node $x^{1}$ in Example 8.8 illustrates this.
Furthermore, it is possible for the radius of ${\cal G}^{\nu}$
to be the eccentricity of only on node in ${\cal G}^{\nu}$;
this occurs when the center consists of only one node.
Can another eccentricity occur for only one node? No, by virtue of 
Conditions 8.1 and 8.2.
There must be at least two nodes for each eccentricity larger than the 
radius.  The proof of this is virtually the same as a proof 
of Lesniak for finite graphs \cite[page 176]{c-l}.

We will need two more results.
They hold except for the trivial case where ${\cal G}^{\nu}$
has only one $(\nu-1)$-section and only one $\nu$-node.

{\bf Lemma 8.10.}  {\em Except for the trivial case just noted, a 
non-boundary bordering $\nu$-node $x^{\nu}$
of a $(\nu-1)$-section  ${\cal S}^{\nu-1}$ has an eccentricity 
that is exactly 
$\omega^{\nu}$ larger than the eccentricity of the internal nodes
of ${\cal S}^{\nu-1}$.}

{\bf Proof.} This follows from Theorem 7.1 
and the fact that any eccentric path 
starting at $x^{\nu}$ and entering ${\cal S}^{\nu-1}$ must pass through
exactly one $(\nu-1)$-tip. $\Box$

An {\em end $(\nu-1)$-section} is a $(\nu-1)$-section having 
exactly one boundary $\nu$-node.

{\bf Lemma 8.11.} {\em Except for the trivial case, 
the eccentricity of the internal nodes of an 
end $(\nu-1)$-section is exactly $\omega^{\nu}$ larger than the 
eccentricity of its boundary $\nu$-node.}

{\bf Proof.}  Any eccentric path of any node of $i({\cal S}^{\nu-1})$
must pass through that boundary $\nu$-node.  
So, the argument of the preceding proof
works again. $\Box$
 
\vspace{.07in}
9. The Center Lies in a $\nu$-Block

This is a known result for finite graphs \cite[Theorem 2.2]{b-h},
\cite[Theorem 2.9]{c-l}, which we now extend transfinitely.
As before, the {\em center} of 
${\cal G}^{\nu}$ is the set of nodes having the minimum eccentricity.
To define a ``$\nu$-block,'' we first define the {\em removal}
of a pristine nonsingleton $\nu$-node $x^{\nu}$ to be the following 
procedure:  $x^{\nu}$ is replaced by two or more singleton $\nu$-nodes,
each containing exactly one of the $(\nu-1)$-tips 
of $x^{\nu}$ and with every 
$(\nu-1)$-tip of $x^{\nu}$ 
being so assigned.  We denote the resulting $\nu$-graph by 
${\cal G}^{\nu}-x^{\nu}$.  Then, a subgraph ${\cal H}$ of ${\cal G}^{\nu}$ 
will be called a 
{\em $\nu$-block} of ${\cal G}^{\nu}$ if ${\cal H}$ is a maximal 
$\nu$-connected subgraph
such that, for every $x^{\nu}$, all the branches of ${\cal H}$ lie 
in the same component of ${\cal G}^{\nu}-x^{\nu}$.  A more explicit way of defining 
a $\nu$-block is as follows:  For any $\nu$-node $x^{\nu}$, 
${\cal G}^{\nu}-x^{\nu}$ consists of one or more components.  Choose one of those components. 
Repeat this for every $\nu$-node, choosing one component for each 
$\nu$-node. Then, take the intersection\footnote{This is the subgraph
induced by those branches, each of which lie in all the chosen 
components.} of all those chosen components.  That intersection may 
be empty, but, if it is not empty, it will be a $\nu$-block 
of ${\cal G}^{\nu}$.
Upon taking all possible intersections of components, one component
from each ${\cal G}^{\nu}-x^{\nu}$, and then choosing the 
nonempty intersections,
we will obtain all the $\nu$-blocks of ${\cal G}^{\nu}$.

Furthermore, we define a {\em cut-node} as a nonsingleton $\nu$-node
$x^{\nu}$ such that ${\cal G}^{\nu}-x^{\nu}$ has two or more components.  
It follows that 
the cut $\nu$-nodes {\em separate} the $\nu$-blocks in the sense that 
any path that terminates at two branches in different $\nu$-blocks must pass 
through at least one cut $\nu$-node.  (Otherwise, the two branches 
would be in the same component of ${\cal G}^{\nu}-x^{\nu}$ for 
every $x^{\nu}$ and therefore
in the same $\nu$-block.)
In summary, we have the following:

{\bf Lemma 9.1.}  {\em The $\nu$-blocks of ${\cal G}^{\nu}$ 
partition ${\cal G}^{\nu}$,
and the cut $\nu$-nodes separate the $\nu$-blocks.}

{\bf Proof.} For each $x^{\nu}$, each branch will be in at least one of 
the components of ${\cal G}^{\nu}-x^{\nu}$, and therefore in 
at least one of the 
$\nu$-blocks.  On the other hand, no branch can be in two different 
$\nu$-blocks because then there would be a 
cut $\nu$-node that separates a branch 
from itself---an absurdity. $\Box$

{\bf Lemma 9.2.}  {\em Each $(\nu-1)$-section ${\cal S}^{\nu-1}$ 
is contained in a 
$\nu$-block.}

{\bf Proof.}  Since every $\nu$-node is pristine, any two branches of 
${\cal S}^{\nu-1}$ are connected through a 
two-ended path of rank no greater than $\nu-1$, and 
that path will not meet any $\nu$-node.  Thus, ${\cal S}^{\nu-1}$ 
will lie entirely within a
single component of ${\cal G}^{\nu}-x^{\nu}$, whatever be the 
choice of $x^{\nu}$.
By the definition of a $\nu$-block, we have the conclusion.  $\Box$

By definition, all the bordering nodes of a $(\nu-1)$-section
${\cal S}^{\nu-1}$ will be  $\nu$-nodes.  Moreover, every $(\nu-1)$-section 
${\cal S}^{\nu-1}$ will have at
least one bordering node, and all the bordering nodes of ${\cal S}^{\nu-1}$
will be nodes of ${\cal S}^{\nu-1}$.  Thus, by Lemma 9.2, every 
$\nu$-block ${\cal H}$
will contain the bordering $\nu$-nodes of its $(\nu-1)$-sections, 
and therefore the rank of ${\cal H}$ is $\nu$.  So, henceforth we 
denote ${\cal H}$ 
by ${\cal H}^{\nu}$.  In general, a $\nu$-node can belong to more than one 
$(\nu-1)$-section and also to more than one $\nu$-block.

{\bf Example 9.3.}  Fig. 6 shows a 1-graph in which the $P_{k}^{0}$
$(k=1,2,3,4,5)$ are endless 0-paths and $v^{1}$, $w^{1}$, $x^{1}$, $y^{1}$,
and $z^{1}$ are 1-nodes.  There are two 1-blocks:  
One of them consists of 
$P_{1}^{0}$ along with $v^{1}$ and $w^{1}$, and the other consists
of the $P_{k}^{0}$ $(k=2,3,4,5)$ along with $w^{1}$, 
$x^{1}$, $y^{1}$, and $z^{1}$.
The only cut 1-node is $w^{1}$.    Also, there are five 0-sections,  
each consisting of one endless 0-path along with its two
bordering 1-nodes.  $\Box$

{\bf Example 9.4.}  The condition that the $\nu$-nodes are pristine is needed 
for Lemma 9.2 to hold.  For example, Fig. 7 shows a 1-graph with 
a nonpristine 1-node $x^{1}$, three 1-blocks,\footnote{Here,
we are extending the definition of a 1-block by requiring that
the elementary tips of $x^{1}$ also be placed in singleton nodes.} 
and two 0-sections.
Branch $b_{1}$ induces one 1-block, branch $b_{2}$ indices another 
1-block, and all the branches of the one-ended 0-path $P^{0}$
induce the third 1-block.  However, $b_{1}$ and $b_{2}$ together
induce a single 0-section ${\cal S}^{0}$, and the branches 
of $P^{0}$ induce another 0-section.  ${\cal S}^{0}$ lies in 
the union of two 1-blocks.  $\Box$

We are finally ready to verify the title of this section concerning 
the {\em center} of ${\cal G}^{\nu}$, which by definition 
is the set of nodes having the minimum eccentricity.  
According to Theorem 8.6, such nodes exist.  Having set up 
appropriate definitions and preliminary results for the 
transfinite case, we can now use a proof that is much the same
as that for finite graphs \cite[Theorem 2.2]{b-h},
\cite[Theorem 2.9]{c-l}.

{\bf Theorem 9.5.}  {\em The center of ${\cal G}^{\nu}$ lies in a $\nu$-block.}

{\bf Proof.}  Suppose the center of ${\cal G}^{\nu}$ lies in two or more 
$\nu$-blocks.  Let ${\cal H}_{1}^{\nu}$ and ${\cal H}_{2}^{\nu}$ 
be two of them.  
By Lemma 9.1, there is a cut $\nu$-node $x^{\nu}$ separating
them.  Let $u$ be an eccentric node for $x^{\nu}$, and let $P_{x,u}$
be an $x^{\nu}$-to-$u$ geodesic.  Thus, $|P_{x,u}|=e(x^{\nu})$.
$P_{x,u}$ cannot contain any node different from $x^{\nu}$
in at least one of ${\cal H}_{1}^{\nu}$ and ${\cal H}_{2}^{\nu}$, say,
${\cal H}_{1}^{\nu}$.  Let $w$ be a center node in ${\cal H}_{1}^{\nu}$
other than $x^{\nu}$, and let 
$P_{w,x}$ be a $w$-to-$x^{\nu}$ geodesic.  Then, $P_{w,x}\cup P_{x,u}$
is a path whose length satisfies $|P_{w,x}\cup P_{x,u}|=
|P_{w,x}|+|P_{x,u}|\geq 1+e(x)$.  This shows that the 
eccentricity of $w$ is greater than the 
minimum eccentricity, that is, $w$ is not a center node---a 
contradiction that proves the theorem.  $\Box$.

\vspace{.07in}
10. The Centers of Cycle-Free $\nu$-Graphs

We now specialize our study to a certain kind of $\nu$-graph
that encompasses the class of transfinite trees as a special case.
(A {\em transfinite $\nu$-tree} is a $\nu$-connected $\nu$-graph
having no loops.)  The kind of $\nu$-graph we now deal with is one 
having no $\nu$-loop that passes through more than one $(\nu-1)$-section.  
All other loops are allowed.  
Let us be more specific.

Because all the $\nu$-nodes are pristine, every loop of rank less than $\nu$ 
must lie within a single $(\nu-1)$-section ${\cal S}^{\nu-1}$, that is,
all its nodes are internal nodes of ${\cal S}^{\nu-1}$.  Such loops 
are allowed.  Moreover, 
a $\nu$-loop might also lie in a single $(\nu-1)$-section ${\cal S}^{\nu-1}$ 
in the sense that all its nodes of ranks less than $\nu$ 
are internal nodes of ${\cal S}^{\nu-1}$ and all its
$\nu$-nodes are bordering nodes of ${\cal S}^{\nu-1}$;  thus, all its 
branches lie in 
${\cal S}^{\nu-1}$.  Such a loop will pass through a closed sequence 
$\{x^{\nu}_{1},x^{\nu}_{2},\ldots,x^{\nu}_{k-1},x^{\nu}_{1}\}$
of bordering ${\nu}$-nodes of ${\cal S}^{\nu-1}$ 
alternating with endless $(\nu-1)$-paths within ${\cal S}^{\nu-1}$.
The possibility of such endless $(\nu-1)$-paths within 
${\cal S}^{\nu-1}$ is implied by Lemma 8.4(b).  Such $\nu$-loops within 
${\cal S}^{\nu-1}$ are also allowed.  On the other hand, it is 
possible in general
for a $\nu$-loop to pass through two or more $(\nu-1)$-sections.
For the sake of a succinct terminology, we shall call the latter
kind of $\nu$-loop a {\em cycle}.  We will henceforth assume that 
${\cal G}^{\nu}$ is so structured that it does not have any cycle and will
say that ${\cal G}^{\nu}$ is {\em cycle-free}.  

In conformity with our definition of an {\em end $(\nu-1)$-section}
as a $(\nu-1)$-section having exactly one boundary $\nu$-node, 
we now define a {\em non-end $(\nu-1)$-section} as a 
$(\nu-1)$-section having two or more boundary $\nu$-nodes.
Note that, when ${\cal G}^{\nu}$ is cycle-free, two $(\nu-1)$-sections
cannot share more than one boundary $\nu$-node, for otherwise 
${\cal G}^{\nu}$ would contain a cycle.  However, still more is implied 
by the cycle-free condition.

{\bf Lemma 10.1.}  {\em Assume that
${\cal G}^{\nu}$ is cycle-free.  Then, ${\cal G}^{\nu}$ has only finitely many 
non-end $(\nu-1)$-sections.}

{\bf Proof.}  If there are no non-end $(\nu-1)$-sections, the conclusion 
is trivially satisfied.  So, assume otherwise, and choose
any non-end $(\nu-1)$-section.  Label it and all its boundary $\nu$-nodes
by ``1.''  Label by ``2'' all the non-end $(\nu-1)$-sections 
that share boundary $\nu$-nodes with that 1-labeled 
$(\nu-1)$-section (if such exist), and label their unlabeled 
boundary $\nu$-nodes by ``2,'' as well.  Each 2-labeled 
section shares exactly one boundary $\nu$-node with the 
1-labeled $(\nu-1)$-section and does not share any 2-labeled $\nu$-node
with any other 2-labeled $(\nu-1)$-section, 
for otherwise ${\cal G}^{\nu}$ would contain
a cycle.  It follows that the number of labeled $\nu$-nodes is no 
less than the number of labeled $(\nu-1)$-sections.  Next, 
label by ``3'' all the non-end $(\nu-1)$-sections that share boundary 
$\nu$-nodes with the 2-labeled sections (if such exist), and label
their unlabeled boundary $\nu$-nodes by ``3,'' as well.  
Again, each 3-labeled $(\nu-1)$-section shares exactly one 
boundary $\nu$-node with exactly one 2-labeled $(\nu-1)$-section and does 
not share any 3-labeled $\nu$-node with any other 
3-labeled $(\nu-1)$-section, for otherwise ${\cal G}^{\nu}$ would 
contain a cycle.  Here, too, it follows that the number of 
labeled $\nu$-nodes is no less than the number of labeled 
$(\nu-1)$-sections.  Continue this way. At each step, 
the number of labeled $\nu$-nodes will be no less
than the number of labeled $(\nu-1)$-sections.  
Since there are only finitely many 
boundary $\nu$-nodes (Condition 8.2) and since these are
the $\nu$-nodes that have been labeled, our conclusion follows.  $\Box$

Our next objective is to replace our cycle-free $\nu$-graph 
${\cal G}^{\nu}$ $(\nu\geq 1)$ by a conventional finite tree (i.e.,
a 0-connected 0-graph having no loop and only finitely many branches), 
a correspondence
that will be exploited in the proof of our final theorem.
The nodal eccentricities for ${\cal T}^{0}$ will be related to the nodal
eccentricities for ${\cal G}^{\nu}$ in a simple way (see Lemma 10.2 below).

In the following, the index $k$ will number the boundary 
$\nu$-nodes of ${\cal G}^{\nu}$ (one number for each) as well 
as sets of non-boundary
bordering $\nu$-nodes.  The index $m$ will number the interiors of non-end
$(\nu-1)$-sections (one number for each interior) as well as sets 
of the interiors of end $(\nu-1)$-sections.

Consider, first of all, a non-end $(\nu-1)$-section ${\cal S}^{\nu-1}$.
Each of its boundary $\nu$-nodes $x^{\nu}_{k}$ is replaced by a 0-node 
$x^{0}_{k}$ having the same index number $k$.  Also, all the 
nonboundary bordering $\nu$-nodes of ${\cal S}^{\nu-1}$ (if such exist)
are replaced by a single 0-node $x_{k'}^{0}$.
Finally, the interior $i({\cal S}^{\nu-1})$ is replaced by 
a single 0-node $y^{0}_{m}$.  A branch is inserted between $y^{0}_{m}$ 
and each of the 
$x_{k}^{0}$ and between $y^{0}_{m}$ and $x_{k'}^{0}$ as well.
Thus, ${\cal S}^{\nu-1}$ is 
replaced by a star 0-graph.  We view $x_{k}^{0}$ (resp. $x_{k'}^{0}$,
resp. $y^{0}_{m}$) as representing $x^{\nu}_{k}$ (resp. every 
nonboundary bordering node $x_{k'}^{\nu}$ of ${\cal S}^{\nu-1}$,
resp. every internal node of ${\cal S}^{\nu-1}$).  
For any $y^{\gamma}\in i({\cal S}^{\nu-1})$, 
we have $d(y^{\gamma},x_{k}^{\nu})=d(y^{\gamma},x_{k'}^{\nu})=
\omega^{\nu}$ (Lemma 8.4(a))
and $d(y_{m}^{0},x_{k}^{0})=d(y_{m}^{0},x_{k'}^{0})=1$.

Next, consider all the end $(\nu-1)$-sections that are incident 
to a single boundary node $x_{k''}^{\nu}$.  
We represent all of their interiors by a single 0-node $y_{m'}^{0}$,
and we insert a branch between $y_{m'}^{0}$ and $x_{k''}^{0}$.
If at least one of those end $(\nu-1)$-sections
has a non-boundary bordering $\nu$-node $x_{k'''}^{\nu}$, we 
represent all of them by another single 0-node $x_{k'''}^{0}$, and we 
insert another branch between $y_{m'}^{0}$ and $x_{k'''}^{0}$.
So, all of these end $(\nu-1)$-sections incident to the chosen 
$x_{k''}^{0}$ are represented either by a single branch 
incident at $x_{k''}^{0}$ and $y_{m'}^{0}$ or by two branches in
series incident at $x_{k''}^{0}$, $y_{m'}^{0}$,
and $x_{k'''}^{0}$.  Here, too, we have replaced the said set of 
end $(\nu-1)$-sections by either a one-branch or two-branch
(elementary) star 0-graph.  Again, we view $x_{k'''}^{0}$ 
(resp. $y_{m'}^{0}$) as representing every 
nonboundary bordering node $x_{k'''}^{\nu}$ 
(resp. every internal node) of the said end $(\nu-1)$-sections.
For any internal node $y^{\gamma}$ in any one of those end 
$(\nu-1)$-sections, we have $d(y^{\gamma},x_{k''}^{\nu})
=d(y^{\gamma},x_{k'''}^{\nu})=\omega^{\nu}$ and 
$d(y_{m'}^{0},x_{k''}^{0})=d(y_{m'}^{0},x_{k'''}^{0})=1$.

We now connect all these star 0-graphs together at their end 0-nodes
in the same way that the $(\nu-1)$-sections are connected together at their 
boundary $\nu$-nodes.  The result is a finite 0-tree ${\cal T}^{0}$. 
Indeed, since 
${\cal G}^{\nu}$ is cycle-free, ${\cal T}^{0}$ has no loops.  Also, 
by Condition 8.2 and Lemma 10.1, ${\cal T}^{0}$ has 
only finitely many branches.  

{\bf Lemma 10.2.}  {\em A node $z$ of any rank in ${\cal G}^{\nu}$ has an 
eccentricity $e(z)=\omega^{\nu}\cdot p$ if and only if its
representative 0-node $z^{0}$ in ${\cal T}^{0}$ has the eccentricity 
$e(z^{0})=p$.}

(Here again, $p$ is a natural number.)

{\bf Proof.}  An eccentric path $P^{\nu}$ of any node $z$ of any rank 
in ${\cal G}^{\nu}$ passes alternately through $(\nu-1)$-sections and 
bordering $\nu$-nodes and terminates at $z$ and an eccentric node 
for $z$.  Because all $\nu$-nodes are pristine, the length 
$|P^{\nu}|$ is obtained by counting the $(\nu-1)$-tips traversed
by $P^{\nu}$ and multiplying by $\omega^{\nu}$ (Lemma 8.3).
Furthermore, corresponding to $P^{\nu}$ there is a unique path 
$Q^{0}$ in ${\cal T}^{0}$ whose nodes $x_{k}^{0}$ and $y_{m}^{0}$
alternate in $Q^{0}$ and represent the bordering nodes $x_{k}^{\nu}$
and interiors of 
$(\nu-1)$-sections ${\cal S}^{\nu-1}_{m}$ traversed by $P^{\nu}$.
Each branch of $Q^{0}$ corresponds to one traversal of a 
$(\nu-1)$-tip in $P^{\nu}$, and conversely.  Thus, we have 
$|P^{\nu}|=\omega^{\nu}\cdot p$ and $|Q^{0}|=p$, where $p$ is the 
number of branches in $Q^{0}$.  Also, since $P^{\nu}$ is an 
eccentric path in ${\cal G}^{\nu}$, $Q^{0}$ is an eccentric path in $Q^{0}$.
Whence our conclusion.  $\Box$

Here is our principal result concerning the centers of cycle-free
$\nu$-graphs.

{\bf Theorem 10.3.}  {\em The center of any cycle-free $\nu$-graph
${\cal G}^{\nu}$ has one of the following forms:
\begin{description}
\item[(a)] A single $\nu$-node $x^{\nu}$.
\item[(b)] The interior $i({\cal S}^{\nu-1})$ of a single $(\nu-1)$-section
${\cal S}^{\nu-1}$.
\item[(c)] The set $i({\cal S}^{\nu-1})\cup\{x^{\nu}\}$, 
where $i({\cal S}^{\nu-1})$ is as in 
(b) and $x^{\nu}$ is one of the bordering $\nu$-nodes of ${\cal S}^{\nu-1}$.
\end{description}
}

{\bf Proof.}  In the trivial case where ${\cal G}^{\nu}$ has just one 
$(\nu-1)$-section and just one $\nu$-node, all the nodes have the
same eccentricity and form (c) holds.  So, consider the case where 
${\cal G}^{\nu}$ has at least two $(\nu-1)$-sections or at least two 
$\nu$-nodes.  Because of Lemmas 8.10 and 8.11, 
neither a non-boundary bordering $\nu$-node
nor the interior of an end $(\nu-1)$-section can be in the
center. Thus, any center node of 
${\cal G}^{\nu}$ is either a boundary $\nu$-node or an internal node  
of a non-end $(\nu-1)$-section.  Also, the correspondence between 
boundary $\nu$-nodes $x_{k}^{\nu}$ and their representative 0-nodes 
$x_{k}^{0}$ in ${\cal T}^{0}$ is a bijection, and so, too, 
is the correspondence
between the interiors $i({\cal S}^{\nu-1}_{m})$ of non-end $(\nu-1)$-sections
and their representative 0-nodes $y_{m}^{0}$ in ${\cal T}^{0}$.
The eccentricities of these entities in ${\cal G}^{\nu}$ are related to
the eccentricities of their 
representatives in ${\cal T}^{0}$ as stated in Lemma 10.2.

We now invoke an established theorem 
for finite 0-trees \cite[Theorem 2.1]{b-h};  
namely, the center of such a tree is either a 
single 0-node of a pair of adjacent 0-nodes.  When the center of 
${\cal T}^{0}$ is a single 0-node $x_{k}^{0}$, form (a) holds.  When
the center of ${\cal T}^{0}$ is a single 0-node $y_{m}^{0}$, form (b)
holds.  Finally, when the center of ${\cal T}^{0}$ is a pair of 
adjacent 0-nodes,
one of them will be a 0-node $x_{k}^{0}$ and the other will be 
a 0-node $y_{m}^{0}$, and thus form (c) holds.  $\Box$

Finally, let us note in passing that, except for the trivial case
mentioned above, the periphery of a cycle-free $\nu$-graph,
defined as consisting of those nodes with the maximum eccentricity,
is comprised either of nonboundary bordering $\nu$-nodes
of end $(\nu-1)$-sections, or the interiors of 
end $(\nu-1)$-sections, or both.

\vspace{.1in}
Figure Captions

\begin{description}

\item{Fig. 1.} A 1-graph consisting 
of a one-way infinite ladder along with two branches, $d$ and $e$, 
connected to infinite extremities of the ladder.  $x_{a}^{1}$ and 
$x_{b}^{1}$ are the only nonsingleton 1-nodes;  all the other
1-nodes are singletons.
\item{Fig. 2.} A 1-loop having two 0-sections.
\item{Fig. 3.} The 1-graph of Example 6.5.
\item{Fig. 4.} The $\nu$-graph ${\cal H}^{\nu}$.  The lines (other than
those of the rectangle) denote isolated endless $(\nu-1)$-paths.
\item{Fig. 5.} The 1-graph of Example 8.9.
\item{Fig. 6.} The 1-graph of Example 9.3.
\item{Fig. 7.} The 1-graph of Example 9.4.
\end{description}

\end{document}